\documentclass[12pt,final]{article}

\usepackage[pdftex,dvips]{graphicx}
\usepackage{epsfig,subfigure} 
\usepackage{dcolumn}      
\usepackage{bm}           
\usepackage{amsmath,amssymb,url,color,colordvi,psfrag,framed,mathtools}

\usepackage{latexsym,amssymb,amsmath,amsbsy,amstext} 
\usepackage{amsfonts}
\usepackage{citesort}
\usepackage{textcomp}
\usepackage{epsfig,subfigure}
\usepackage{multirow}
\usepackage{array}
\usepackage{wrapfig}
\usepackage{psfrag,color,colordvi}
\usepackage{bbm}

\newcommand{\eref}[1]{(\ref{#1})}
\newcommand{\fref}[1]{Figure~\ref{#1}}
\newcommand{\tref}[1]{Table~\ref{#1}}
\newcommand{\sref}[1]{Section~\ref{#1}}
\newcommand{\vm}[1]{\mathbf{#1}}
\newcommand{\vx}{\vm{x}}

\newcommand{\R}{\mathbb{R}}

\newcommand{\rhop}{\rho^+}
\newcommand{\rhom}{\rho^-}
\newcommand{\rhot}{\tilde\rho}
\newcommand{\rhomt}{\tilde\rho^-}
\newcommand{\rhopt}{\tilde\rho^+}
\newcommand{\tauv}{\boldsymbol{\tau}}
\newcommand{\Vp}{V^+}

\newcommand{\vImt}{\tilde{v}_I^-}
\newcommand{\Vmt}{\tilde{V}^-}
\newcommand{\Vpt}{\tilde{V}^+}
\newcommand{\valmt}{\tilde{v}_\alpha^-}
\newcommand{\Valmt}{\tilde{V}_\alpha^-}
\newcommand{\rhoIm}{\rho_I^-}
\newcommand{\rhoImt}{\tilde\rho_I^-}
\newcommand{\rhoIt}{\tilde\rho_I}
\newcommand{\Rv}{\mathbf{R}}
\newcommand{\mt}{M{\"o}bius transformation~}

\begin{document}
\title{Efficient adaptive integration of functions with sharp gradients and
cusps in $n$-dimensional parallelepipeds}

\author{S. E. Mousavi$^1$, J. E. Pask$^2$ and 
N. Sukumar$^{1,}$\footnote{Corresponding author. 
E-mail: \texttt{nsukumar@ucdavis.edu}}}

\date{
$\mbox{}^1$Department of Civil and Environmental Engineering \\
University of California, Davis, CA 95616. \\
$\mbox{}^{2}$Condensed Matter and Materials Division\\
Lawrence Livermore National Laboratory, Livermore, CA 94550.\\
\vspace*{0.2in}
\today}

\maketitle

\begin{abstract}
In this paper, we study the efficient numerical integration of functions with
sharp gradients and cusps. An adaptive integration algorithm is presented that
systematically improves the accuracy of the integration of a set of functions.
The algorithm is based on a divide and conquer strategy and is independent of
the location of the sharp gradient or cusp. The error analysis reveals that for
a $C^0$ function (derivative-discontinuity at a point), a rate of convergence of
$n+1$ is obtained in $\R^n$. Two applications of the adaptive integration scheme
are studied. First, we use the adaptive quadratures for the integration of the
regularized Heaviside function---a strongly localized function that is used
for modeling sharp gradients. Then, the adaptive quadratures are employed in
the enriched finite element solution of the all-electron Coulomb problem in
crystalline diamond. The source term and enrichment functions of
this problem have sharp gradients and cusps at the nuclei. We show that the
optimal rate of convergence is obtained with only a marginal increase in the
number of integration points with respect to the pure finite element solution
with the same number of elements. The adaptive integration scheme is
simple, robust, and directly
applicable to any generalized finite element method employing enrichments with
sharp local variations or cusps in $n$-dimensional parallelepiped elements.
\end{abstract}

\section{Introduction}
Functions with sharp gradients appear in the solution of problems with
localization and cohesive process
zones~\cite{patzak:2003:PZR,benvenuti:2008:ARX,benvenuti:2008b:ARX},
shear bands~\cite{areias:2006:TSS,areias:2007:TSM}, 
thermal
gradients~\cite{tamma:1989:HPF,merle:2002:STA,ohara:2009:GFE}, 
convection- and advection-diffusion
problems~\cite{brooks:1982:SUP,abbas:2010:TXF,kalashnikova:2011:ADE}, and in
electronic structure
calculations and {\em ab initio} materials
modeling~\cite{challacombe:2000:LSC,sukumar:2009:CAE,pask:2011:LSS}.
In the partition-of-unity~\cite{melenk:1996:PUM} solution of these problems,
such sharp gradients and cusps are efficiently resolved by
incorporating enrichment functions
that resemble the solution locally. As a result, efficient numerical
integration of the basis functions and their gradients to form the system
matrices becomes computationally demanding since one has to deal with strongly
localized functions, instead of polynomial integrands. In many applications,
the enrichment functions are the solution of local problems and known only
numerically. Evaluation of such integrands can be extremely time-consuming,
which points to the need for an efficient integration scheme. In this paper,
we present an adaptive scheme for the integration of functions with sharp
 gradients and cusps. Adaptive integration algorithms accumulate integration 
points in
regions with higher errors, and use fewer points where the integrand is smooth.
Also, we are interested in domains that can be 
 prescribed as a collection of
hyperparallelepipeds, for example, parallelograms and parallelepipeds in two
and three dimensions, respectively.

Adaptive integration schemes are normally recursive in nature and have a few
common ingredients~\cite{berntsen:1993:AAF}: a quadrature rule that can be
applied to the integration domain to provide a local estimate of the
integration; a procedure to estimate the local integration error; a strategy
to partition the integration domain into smaller divisions of the same shape;
and a stopping criterion. The algorithm of
Gander and Gautschi~\cite{gander:2000:AQR} over the interval and that of
Berntsen et al.~\cite{berntsen:1992:AAF,berntsen:1993:AAF} for a collection
of triangles and tetrahedra are examples.
Genz and Cools~\cite{genz:2003:AAN} proposed an algorithm for a vector-valued
function over a combination of $n$-dimensional simplices. At each step, a subset
of the simplices with the highest errors are selected and subdivided, and
quadratures over the subdivisions are used to update the integral. The local
error in Reference~\cite{genz:2003:AAN} is estimated by applying
\emph{null} quadrature rules---quadratures that
(incorrectly) integrate to zero all polynomials up to degree $d$, and fail to
do so for at least one polynomial of a higher degree
$d+1$~\cite{lyness:1965:SIR,berntsen:1991:EEI}. Herein, we
use the difference in the
integration of a function with two different tensor-product quadrature rules
as the local error estimate, and proceed to subdivide a cell until the absolute
error of integration of all the functions over each individual cell falls below
a prescribed tolerance.

In this paper, we consider two types of local features. First, we focus on
the regularized Heaviside function: this function has a sharp gradient,
and by shrinking the size of the process zone, becomes strongly localized.
Tornberg~\cite{tornberg:2002:MDQ} and
Patz{\'a}k and Jir{\'a}sek~\cite{patzak:2003:PZR} proposed regularized forms
of the Heaviside function that smear the strong discontinuity over a short
distance on which the discontinuous function is approximated by a localized
function. Oh et al.~\cite{oh:2007:TSP} introduced several
smooth-piecewise-polynomial regularized discontinuous functions that can be
used in the extended finite element method.
Benvenuti et al.~\cite{benvenuti:2008:ARX} presented a method to integrate the
regularized Heaviside function by the integration of the equivalent smooth
functions. This is an extension of the method of Ventura~\cite{ventura:2006:OEQ}
 for the integration of discontinuous functions, which is however restricted to
 elements with constant Jacobian of the transformation. A more general
 treatment for the integration of arbitrary classes of functions, which is
 also adopted in this work, is adaptive integration using a posteriori error
 estimates: the integration points are concentrated close to the region where
 sharp gradients appear, and fewer points are used
 elsewhere~\cite{benvenuti:2008:ARX}.

Next, we consider functions with cusps. We pick the Coulomb problem in
crystalline diamond and apply an enriched finite element (EFE) approach: Poisson's
equation is solved with the electronic charge density as the source
term, and the isolated atom solutions as the enrichment 
functions~\cite{pask:2011:LSS}.
These functions have sharp gradients in the region close to the atomic
sites, and have cusps at the nuclei. In one dimension, \mt has been
used for the integration of functions with a peak at or near a
boundary~\cite{homeier:1990:NIO,homeier:1992:OTE,lopez:1994:COT}.
When \mt is applied to a standard quadrature over the interval, integration
points are attracted toward the sharp gradient and more accurate results are
obtained. However, this technique cannot be applied in higher dimensions. A
well-known remedy in such cases is adaptive integration. Our adaptive
integration scheme recursively performs a uniform refinement of the
parent cell until the prescribed error tolerance over each subcell
is met. The number of integration points in the adaptive quadratures
is only marginally more than that obtained for the finite element solution
with the same number of elements (and much lower accuracy), which points to
the fact that numerical integration is not a bottleneck for the EFE solution
of the problem.

The structure of this paper follows. In~\sref{sec:algorithm}, the adaptive
integration scheme is introduced, followed by an error analysis
in~\sref{sec:error}. It is shown that the algorithm is efficient for
integrands with a cusp, such as $C^0$ functions of the form
$f(\vx) = 1-r$ or $f(\vx) = \exp(-\alpha r)$. The
adaptive integration algorithm is used for the integration of the regularized
Heaviside function in~\sref{sec:heaviside}. In~\sref{sec:poisson}, an enriched
finite element approach is applied to solve the all-electron 
Coulomb problem in crystalline diamond 
(charge density has a cusp at the nuclei).  We close
with a few concluding remarks in~\sref{sec:conclusions}.

\section{Adaptive Integration Scheme}\label{sec:algorithm}
Our quadrature construction algorithm is customized to meet the integration
demands of high accuracy over parallelepipeds. A tensor-product quadrature
based
on a one-dimensional Gauss rule with five points in each direction is used to
evaluate the local integrals. A tensor-product quadrature with eight points
in each direction is used to compute the reference integral, and provides an
estimate of the local
integration error. The reason for choosing eight points in each direction is
that a quadrature with six points might not be accurate enough to detect the
integration error (due to the sharp gradients of the integrands);
and a seven-point quadrature may suffer from the same
odd-even defect of the quadratures as does a five-point quadrature.
If the absolute error of integration is greater than the
prescribed tolerance, the integration domain is uniformly divided into eight
cells (in three dimensions) and the adaptive integration is performed over
each cell recursively.
This process is started with all the functions that need to be integrated,
and at each step only those functions whose integration error is greater
than the tolerance are passed to the next level, until all functions are
integrated to the specified tolerance. Using the relative error as the
stopping criterion may cause numerical difficulties in parts of the domain
where the integral of at least one of the functions is close to zero. This
problem was also reported in Reference~\cite{gander:2000:AQR}, and was
resolved by additionally checking the local estimate against the global
estimate of the integration: if the contribution of the partition was small
enough, further subdivision was circumvented. Such a problem is avoided by using
the absolute error as the measure of accuracy. Note that the error tolerance
in our quadrature construction scheme is only a stopping criterion, and may
not represent the exact integration error. However, it provides a systematic
means to improve the accuracy of the numerical integration until stable
solutions are obtained.  
In a specific application, the appropriate tolerance depends on the 
overall accuracy required and the mesh resolution: as the mesh is refined, the
tolerance (error per element) should be decreased proportionally to attain the
same overall accuracy. Furthermore, if material parameters that appear in the
weak form integrals have strong spatial variations, they can be included 
in the integrands $f_i$ (see Algorithm below) that are evaluated to construct 
the quadrature rule within the element.
The following pseudo-code explains our quadrature
construction scheme.

\noindent\textbf{Algorithm: Adaptive quadrature construction} \\
{\bf Input:} The domain of integration $\Omega$; $numf$ integrands $\{f_i\}$;
prescribed error tolerance $tol$ \\
{\bf Output:} Adaptive quadrature over $\Omega$ that integrates
$\{f_i\}_{i=1}^{numf}$ within the accuracy of $tol$ \\
for $i = 1:numf$ do \\
\indent\indent $I_1$ =  integrate $f_i$ over $\Omega$ with $5$-point quadrature \\
\indent\indent $I_2$ =  integrate $f_i$ over $\Omega$ with $8$-point quadrature \\
\indent\indent if $|I_2-I_1| > tol$ then set $P_i = 1$, otherwise set $P_i = 0$ \\
end do \\
if $\{P_i\}_{i=1}^{numf} = 0$ \\
\indent\indent return the $5$-point quadrature as the rule over $\Omega$ \\
else \\
\indent\indent partition $\Omega$ into eight uniform cells $\{\Omega_j\}_{j=1}^8$ \\
\indent\indent for each cell $\Omega_j$, call the quadrature construction routine with the integrands \\
\indent\indent $\{f_i$, s.t. $P_i = 1\}_{i=1}^{numf}$ \\
\indent\indent put the eight obtained quadratures together: the adaptive
               quadrature over $\Omega$ \\
end if

This adaptive integration scheme is similar to that of
van Dooren and de Ridder~\cite{dooren:1976:AAA}, except that in
Reference~\cite{dooren:1976:AAA}, the domain is subdivided in one
direction (dimension) at a time, and all integrands are carried until
the last iteration. Limiting subdivision to one dimension at each step
can be useful if the integrands have a strong dependence on one of the
spatial dimensions. Berntsen et al.~\cite{berntsen:1991:AAA} devise an
analogous scheme with nonuniform subdivision, and improve the error
estimate: the error is approximated by combining integrals over a cell
and its children. Adaptive integration by octree subdivision is also used
by Pieper~\cite{pieper:1999:RGI}, and the integration error is estimated
as the difference between the integral over a cell and its children.

\subsection{Example}\label{sec:adap_example}
\begin{figure}
\setlength{\unitlength}{0.0625\linewidth}
\centering
\begin{picture}(16,4)
  \put(.1,.2){\includegraphics[width=3\unitlength]{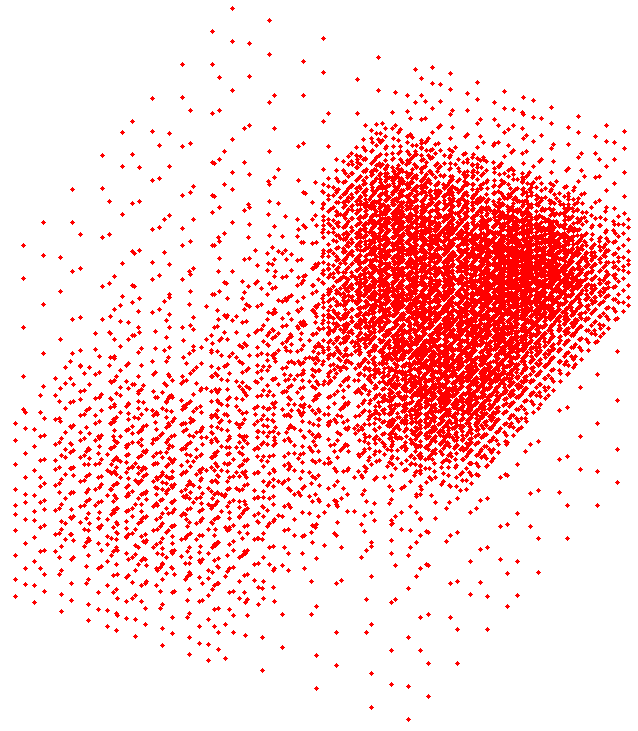}}
  \put(3.3,.2){\includegraphics[width=3\unitlength]{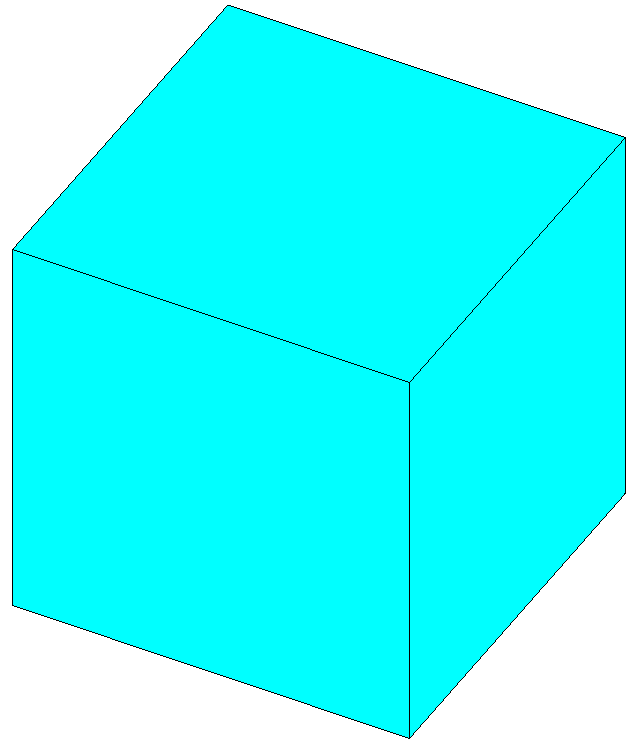}}
  \put(6.5,.2){\includegraphics[width=3\unitlength]{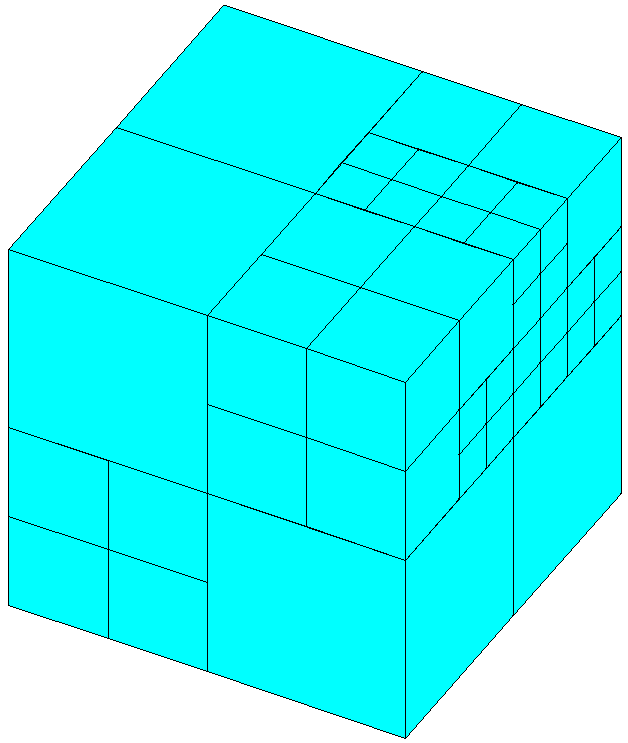}}
  \put(9.7,.2){\includegraphics[width=3\unitlength]{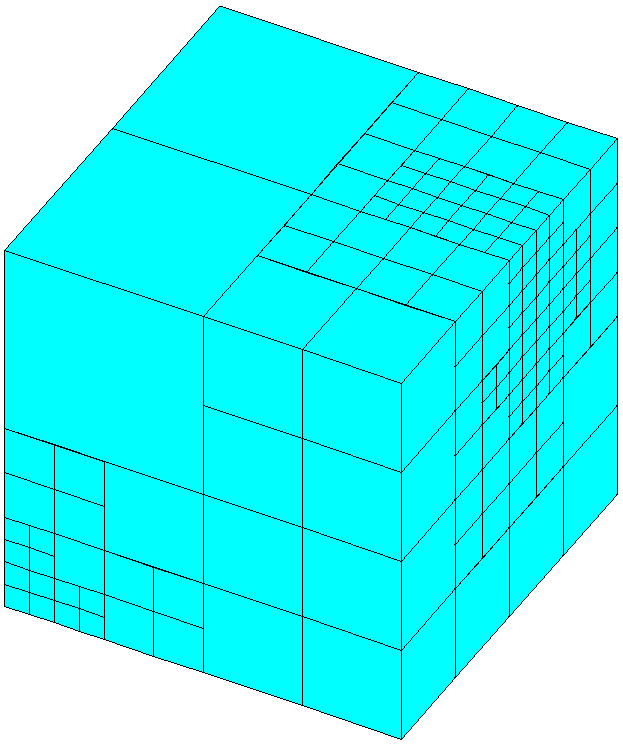}}
  \put(12.9,.2){\includegraphics[width=3\unitlength]{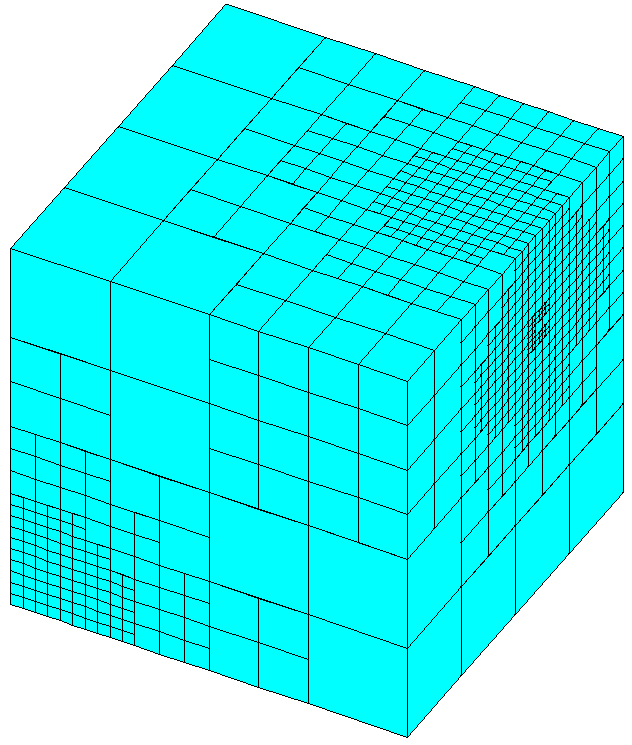}}
  \put(1.6,.2){\makebox(0,0)[tc]{(a)}}
  \put(4.8,.2){\makebox(0,0)[tc]{(b)}}
  \put(8,.2){\makebox(0,0)[tc]{(c)}}
  \put(11.2,.2){\makebox(0,0)[tc]{(d)}}
  \put(14.4,.2){\makebox(0,0)[tc]{(e)}}
\end{picture}
\caption{Adaptive quadrature over the unit cube.
(a) quadrature points; and
(b)-(e) the sequence of octree refinements.}\label{fig:adaptive_example}
\end{figure}
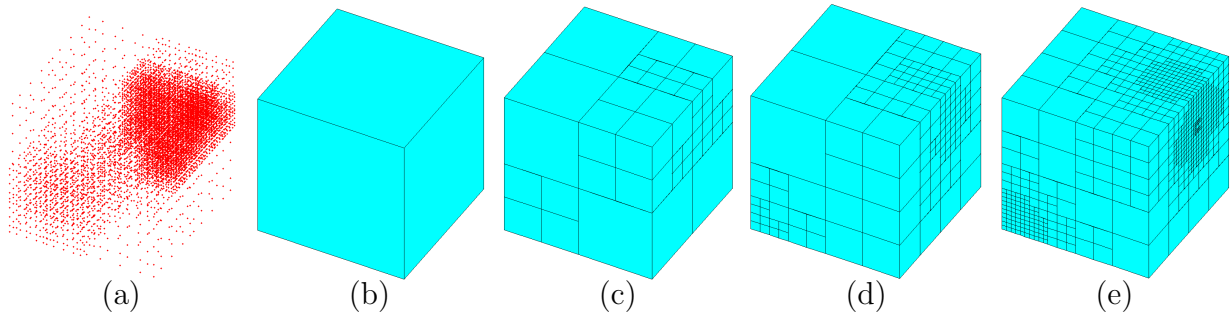

Consider the functions $f_1(\vx) = 10\exp(-100r_1^2)$ and
$f_2(\vx) = 100\exp(-200r_2^2)$ as the integrands for the quadrature
construction, where $r_1$ is the distance from the point $\vx$ to the
origin, $r_2$ is the distance to $(.81, .62, .73)$, and the domain of
integration is the unit cube. Figure~\ref{fig:adaptive_example}a shows
an adaptive quadrature over the unit cube with the above integrands and
an error tolerance of $10^{-6}$ with 8875 integration points. A sequence
of octree refinements of the domain is shown
in~\fref{fig:adaptive_example}b--\ref{fig:adaptive_example}e.
The \texttt{MATLAB\texttrademark} implementation of the adaptive quadrature
is listed in Appendix~\ref{sec:matlab_code}.

\section{Error Analysis}\label{sec:error}
In this section, the convergence properties of the adaptive integration
algorithm are studied. Specifically, we are interested in cases where the
integrand has a cusp, for example $f(\vx) = 1-r$ or
$f(\vx) = \exp(-\alpha r)$, which are $C^0$
functions with a derivative-discontinuity at the origin. However, this will not
have an adverse effect on the efficiency of the adaptive integration scheme,
and a high rate of convergence is realized. The following examples clarify
the problem.

\begin{figure}
\setlength{\unitlength}{0.0625\linewidth}
\centering
\begin{picture}(16,6.5)
  \put(0,.2){\includegraphics[width=7\unitlength]{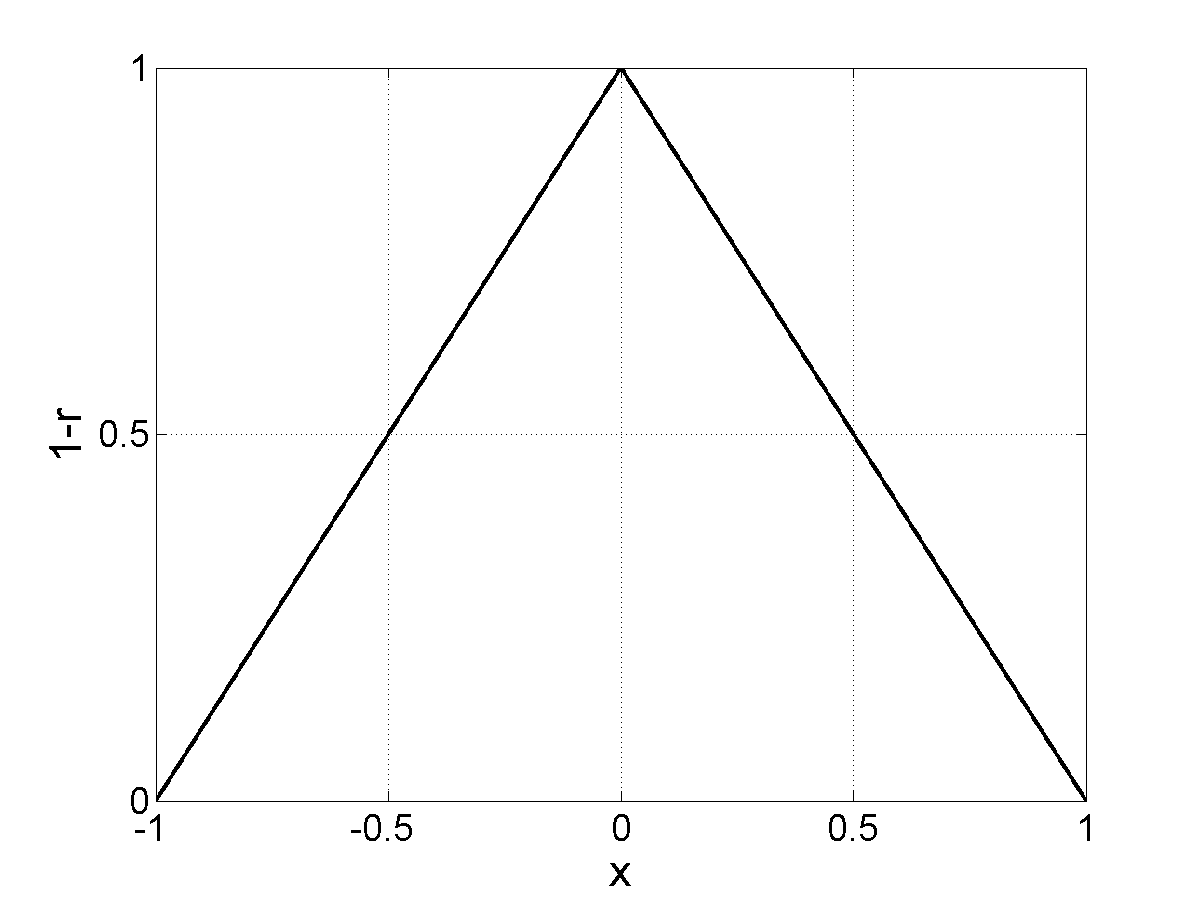}}
  \put(7,0){\includegraphics[width=8.5\unitlength]{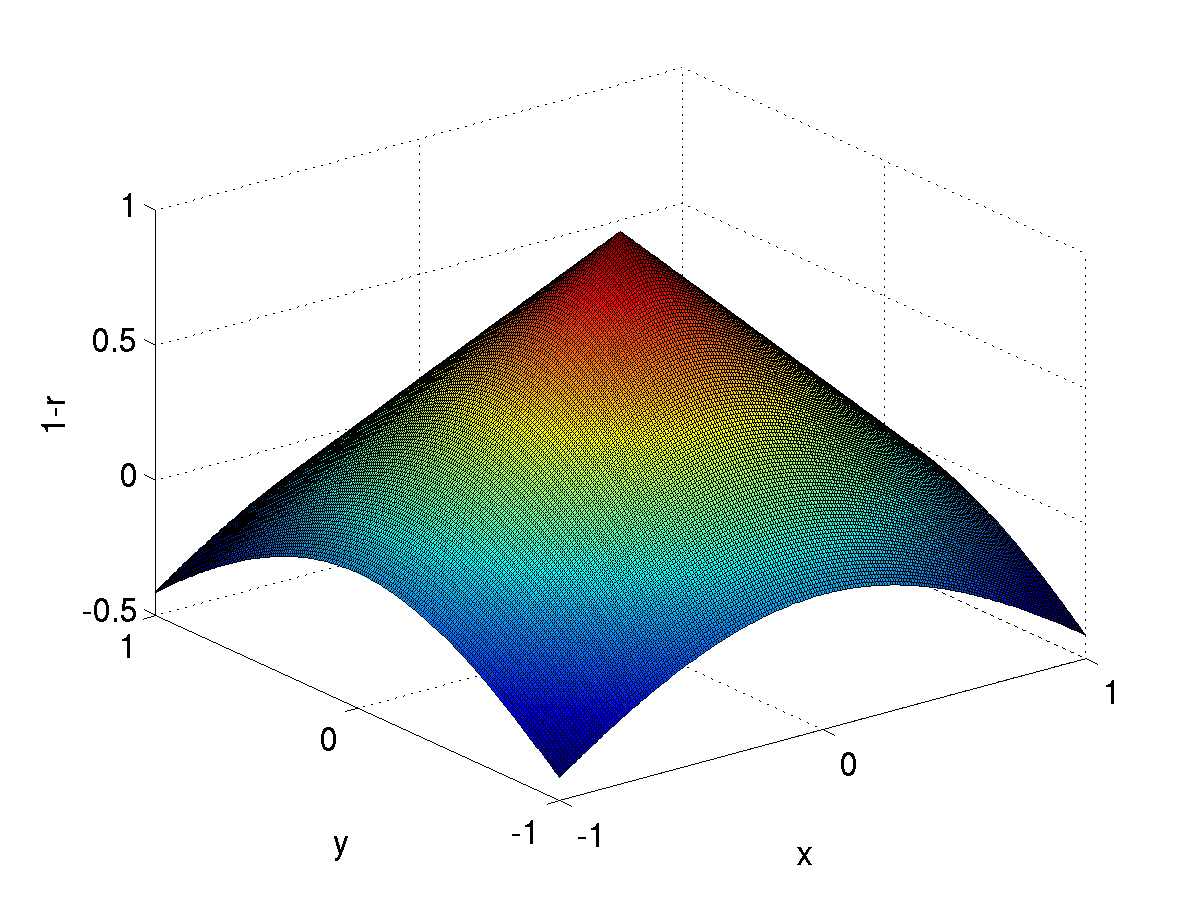}}
  \put(3.5,0){\makebox(0,0)[tc]{(a)}}
  \put(11.5,0){\makebox(0,0)[tc]{(b)}}
\end{picture}
\caption{$f(\vx) = 1-r$.
(a) in one dimension; and
(b) in two dimensions.}\label{fig:adaptive_cusps}
\end{figure}

\begin{figure}
\setlength{\unitlength}{0.0625\linewidth}
\centering
\begin{picture}(16,6)
  \put(0,.2){\includegraphics[width=8\unitlength]{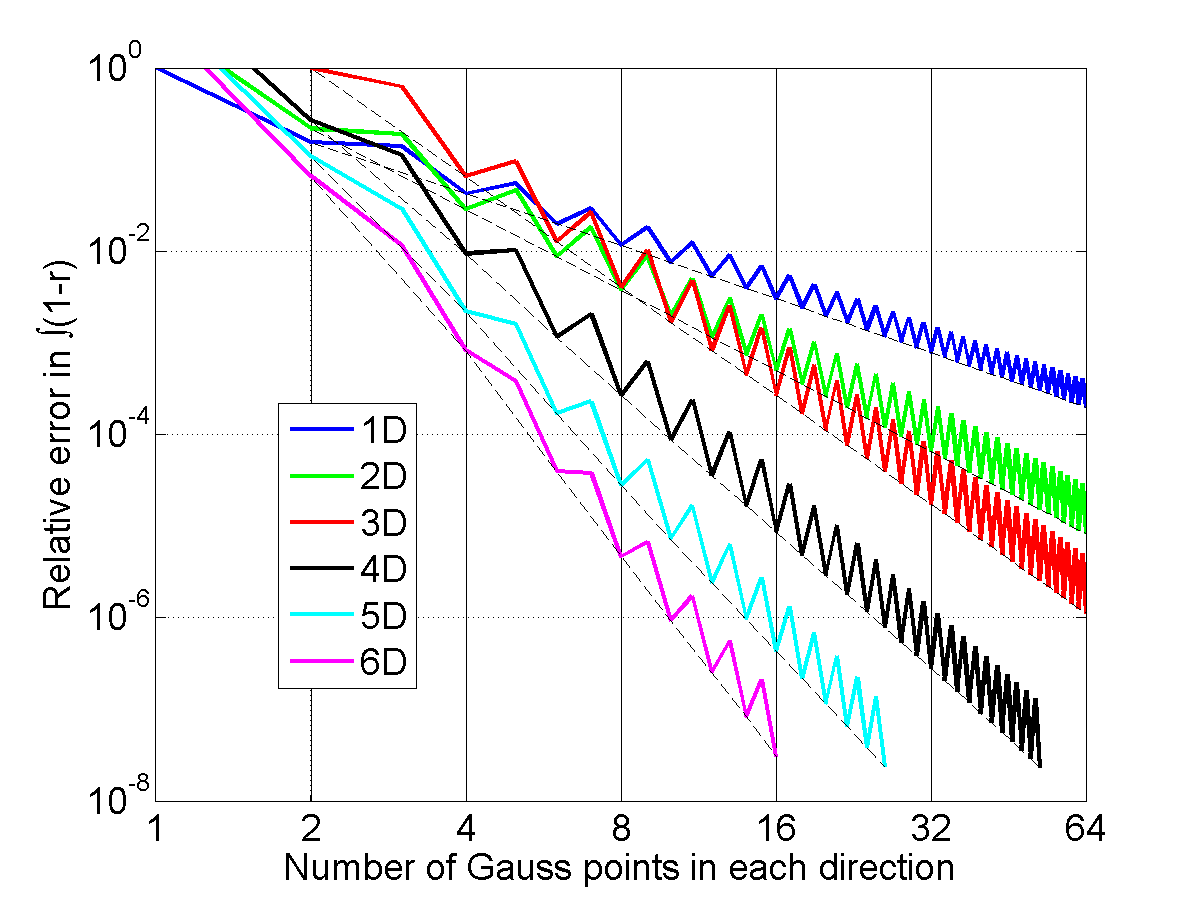}}
  \put(8,.2){\includegraphics[width=8\unitlength]{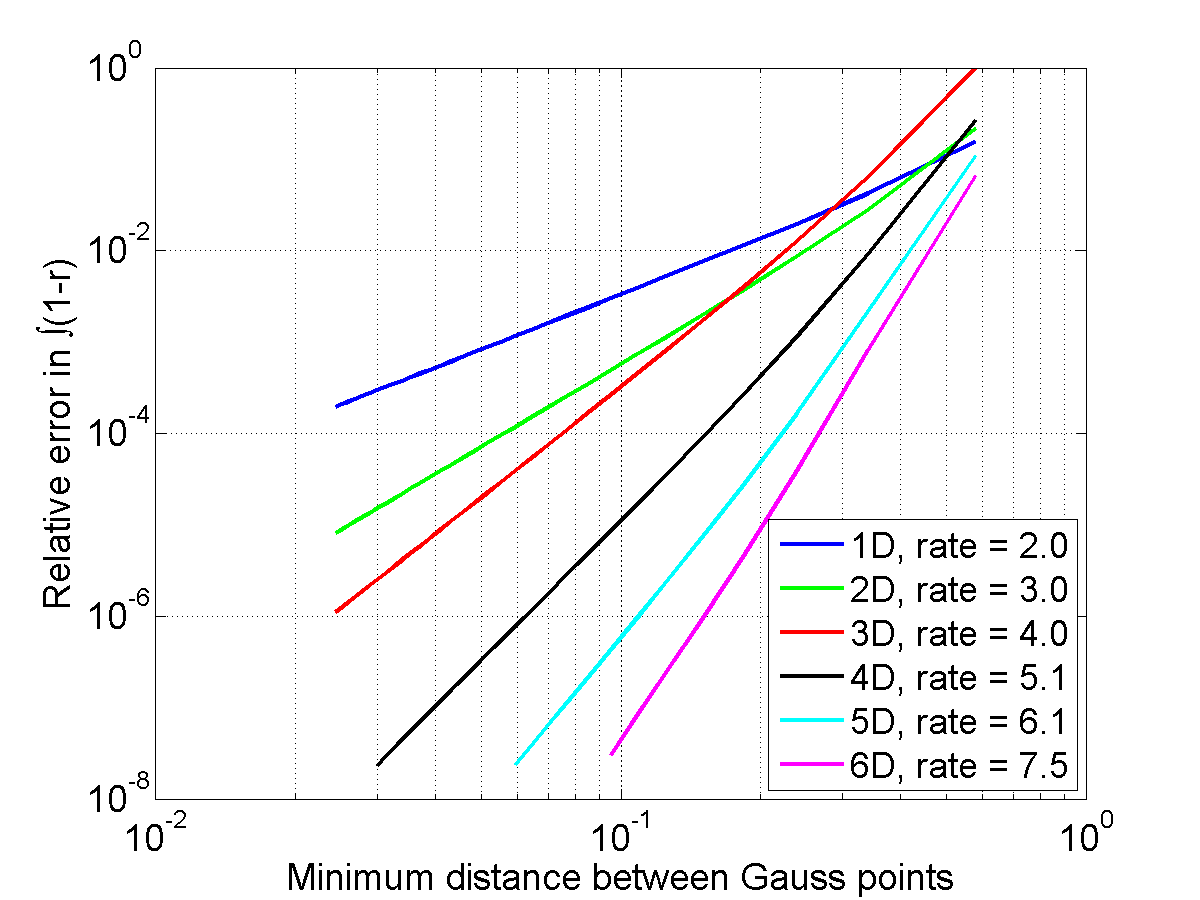}}
  \put(4,0){\makebox(0,0)[tc]{(a)}}
  \put(12,0){\makebox(0,0)[tc]{(b)}}
\end{picture}

\vspace*{0.1in}
\caption{Convergence curves of the integration of $f(\vx) = 1-r$.
(a) with respect to the number of integration points in each direction; and
(b) with respect to the minimum distance of the Gauss points to the cusp (only
even number of integration points are shown).}\label{fig:adaptive_cusps_convergence}
\end{figure}

\begin{figure}
\setlength{\unitlength}{0.0625\linewidth}
\centering
\begin{picture}(16,6.5)
  \put(0,.2){\includegraphics[width=8\unitlength]{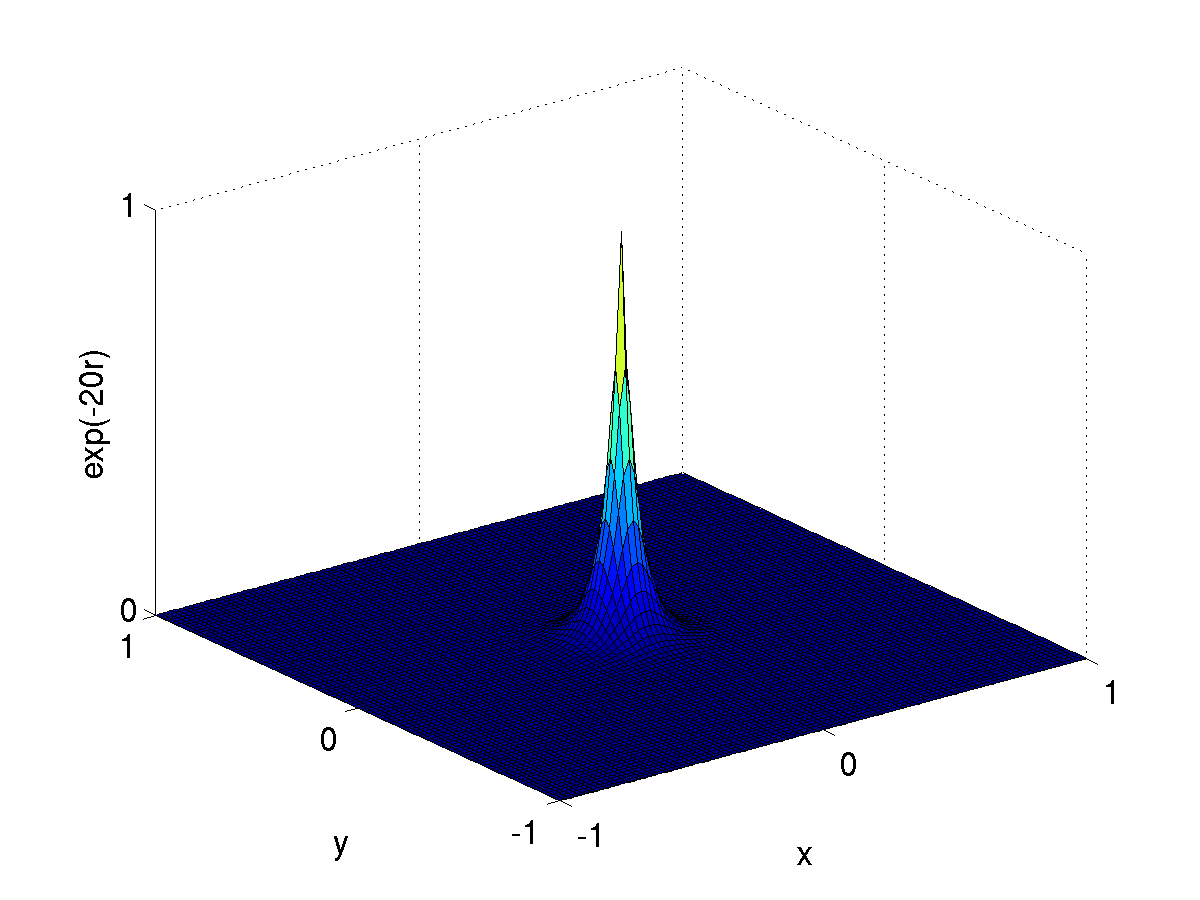}}
  \put(8,.2){\includegraphics[width=8\unitlength]{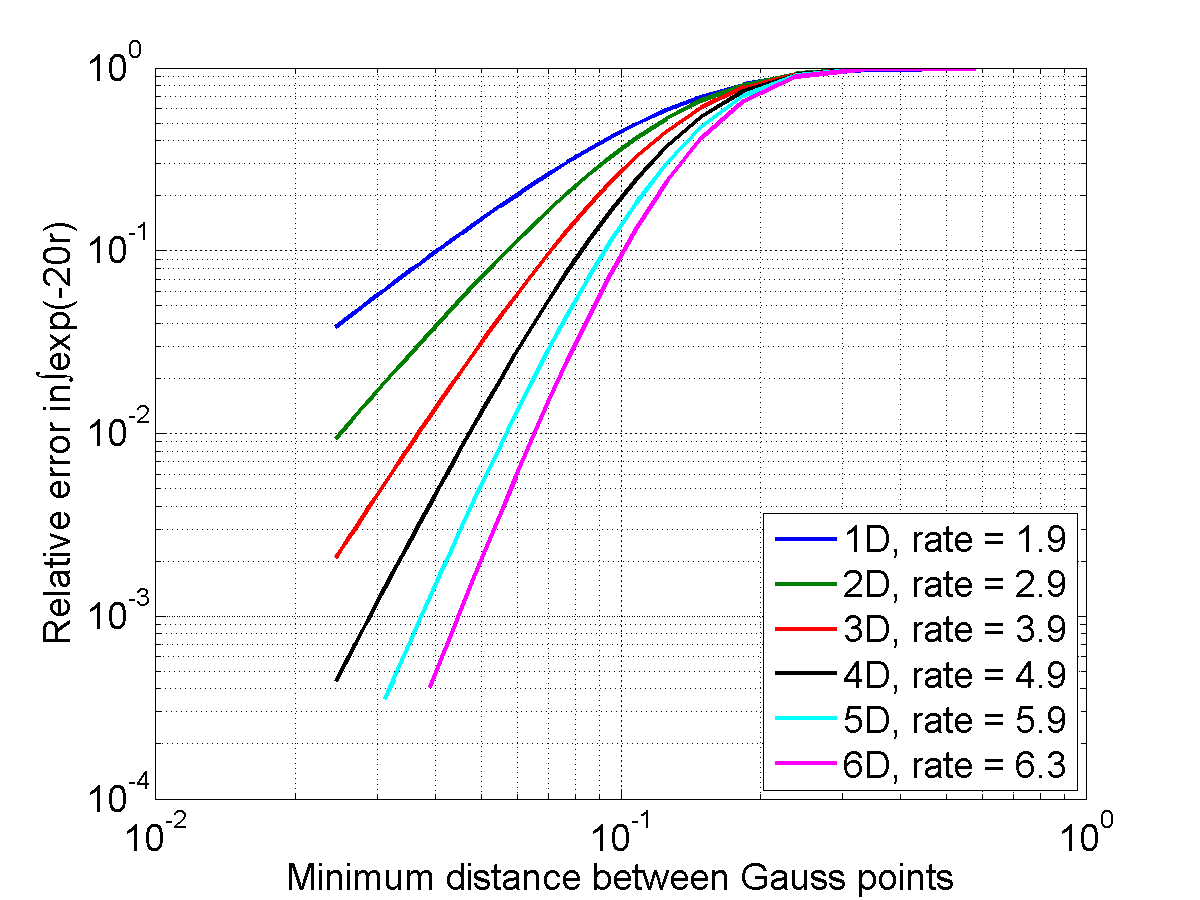}}
  \put(4,0){\makebox(0,0)[tc]{(a)}}
  \put(12,0){\makebox(0,0)[tc]{(b)}}
\end{picture}

\vspace*{0.1in}
\caption{$f(\vx) = \exp(-20r)$.
(a) $f$ in two dimensions; and (b) convergence curves of the integration of
$f(\vx)$ with respect to the minimum distance of the Gauss points to the cusp
(only even number of integration points are shown).}\label{fig:adaptive_cusps_exp}
\end{figure}

Consider the integration of the function $f(\vx) = 1-r$,
where $r$ is the distance of the point $\vx$ to the origin. The
integration domain is $[-1,1]^n$ for $n = 1, 2, \ldots, 6$. 
In~\fref{fig:adaptive_cusps}, the function $f$
is illustrated in one and two dimensions. Gauss
quadratures are used to integrate $f$, and the position of the cusp is arbitrary.
By contrast, one could consider the cusp as a hindrance to convergence, and
partition the integration domain by placing the cusp at the boundary of the
subdivisions. \fref{fig:adaptive_cusps_convergence} shows the convergence of
the integration as the number of integration points is increased: a rate of
convergence of $n+1$ is observed in $\R^n$,
which can be explained by appealing to the Taylor
expansion of the integrand.
Since the integrand is continuous over a small
cell $\Omega_0$ containing the cusp, the error of approximating it with a
polynomial is $\mathcal{O}(h)$; hence the error in the integration is
$\mathcal{O}(h) V_{\Omega_0} \sim \mathcal{O}(h^{n+1})$, where
$V_{\Omega_0} \sim \mathcal{O}(h^n)$ is the volume of $\Omega_0$. This
indicates that our quadrature construction scheme can
attain higher convergence rates in higher dimensions (for $C^0$
functions), in contrast to the usual deficiency of numerical integration
methods that suffer from the curse of dimensionality.
We perform the same experiment with $f(\vx) = \exp(-20r)$, which has a cusp at
the origin. The function is plotted in~\fref{fig:adaptive_cusps_exp}a in two
dimensions. The convergence curves for the integration using tensor-product
quadratures are shown in~\fref{fig:adaptive_cusps_exp}b with respect to the
minimum distance of the integration points to the cusp. As expected, a rate of
convergence of $n+1$ is observed in $\R^n$.
Numerical experiments also reveal that for smoother functions $C^m$, $m > 0$,
higher convergence rates are achieved, consistent with
Darboux's Principle~\cite{boyd:2001:CAF}.


\section{Numerical Examples}\label{sec:examples}
\subsection{Integration of the Regularized Heaviside Function}\label{sec:heaviside}
\begin{figure}
\setlength{\unitlength}{0.0625\linewidth}
\centering
\begin{picture}(16,10.5)
  \put(2,7.2){\includegraphics[width=3\unitlength]{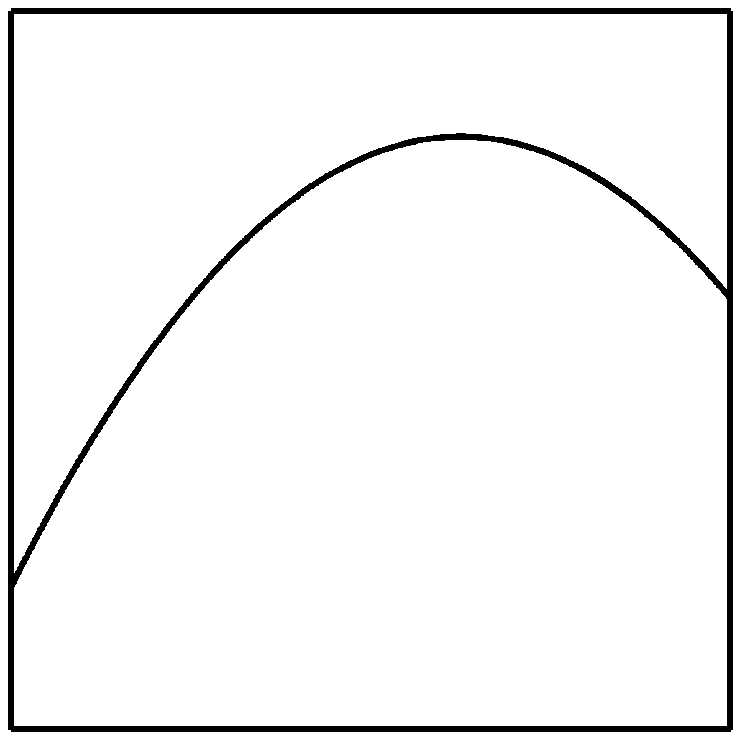}}
  \put(6.5,7.2){\includegraphics[width=3\unitlength]{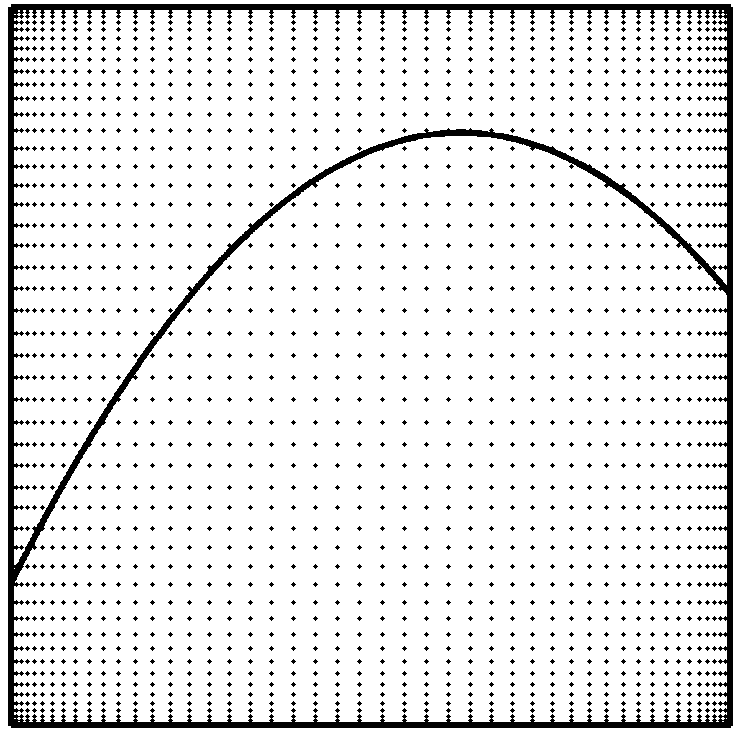}}
  \put(11,7.2){\includegraphics[width=3\unitlength]{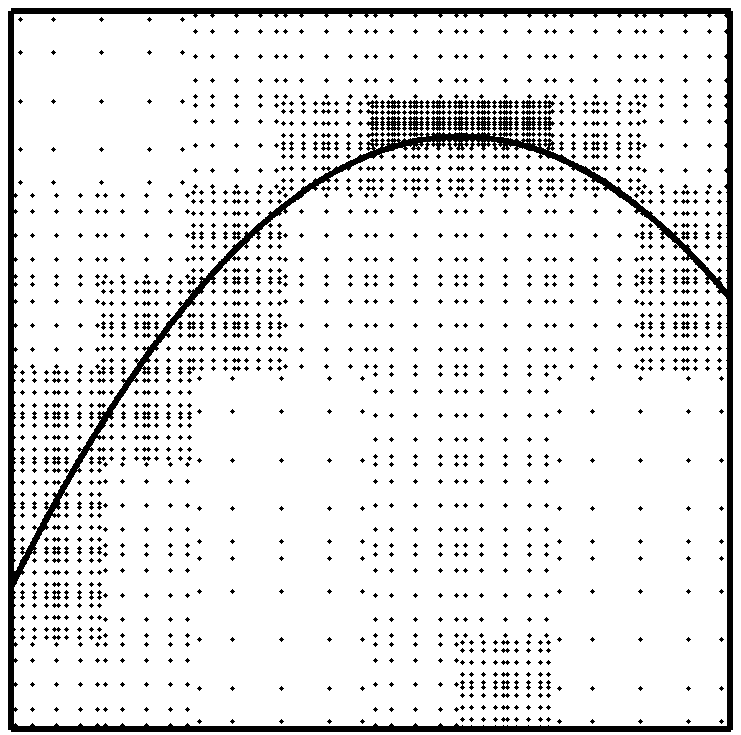}}
  \put(0,.2){\includegraphics[width=8\unitlength]{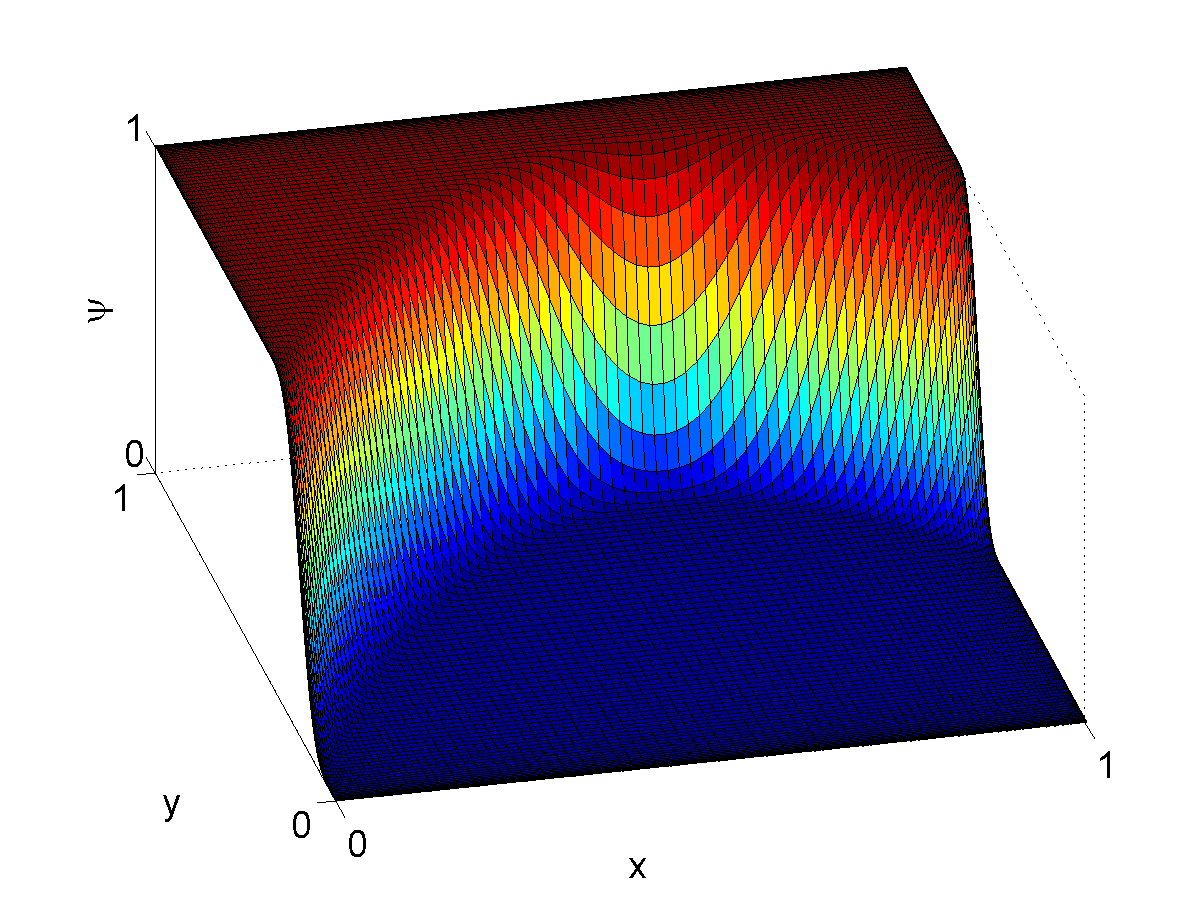}}
  \put(8,.2){\includegraphics[width=8\unitlength]{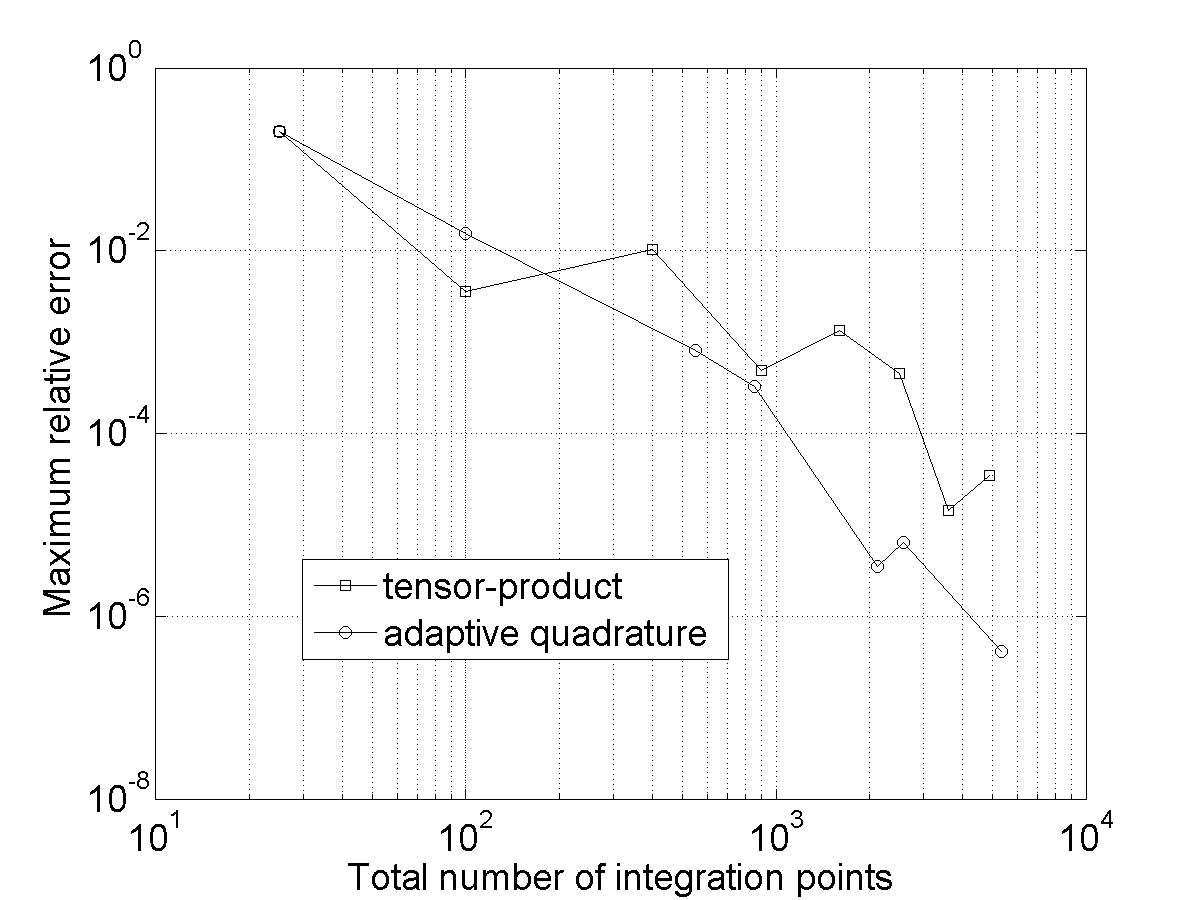}}
  \put(3.5,7){\makebox(0,0)[tc]{(a)}}
  \put(8,7){\makebox(0,0)[tc]{(b)}}
  \put(12.5,7){\makebox(0,0)[tc]{(c)}}
  \put(4,0){\makebox(0,0)[tc]{(d)}}
  \put(12,0){\makebox(0,0)[tc]{(e)}}
\end{picture}

\vspace*{0.1in}
\caption{Numerical integration of the regularized Heaviside function with a
curved interface.
(a) domain of integration and interface; (b) tensor-product;
(e) adaptive quadrature;
(f) the regularized Heaviside function for $\varepsilon = 0.085$; and
(g) maximum relative error in the integration of the five functions versus
the number of integration points for different strategies.}\label{fig:highgradient_ex3}
\end{figure}

Solution fields with sharp spatial gradients arise in modeling physical
phenomena such as shear band evolution, damage and cohesive process zones,
and convection-dominated problems with shocks. Different classes of the
regularized step function are adopted as the enrichment function for
modeling sharp gradients, among which the following piecewise polynomial
regularized Heaviside function is proposed by
Patz{\'a}k and Jir{\'a}sek~\cite{patzak:2003:PZR}:
\begin{equation} \label{eq:regularized_heaviside}
\psi(\phi, \varepsilon) = \left\{\begin{array}{ll} 0 & \mbox{if } \phi < -\varepsilon \\
\frac{1}{V_{\varepsilon}} \int_{-\varepsilon}^{\phi} \left(1 - \frac{\xi^2}{\varepsilon^2}\right)^4 d\xi & \mbox{if } |\phi| \le \varepsilon \\
1 & \mbox{if } \phi > \varepsilon, \end{array}\right.
\end{equation}
where $\phi(\vx)$ is the signed distance from the interface, and $\varepsilon$
is a parameter that determines the gradient of the enrichment function
(half-width of the zone). The reference volume $V_{\varepsilon}$ is set to
$256\varepsilon / 315$ to enforce $C^4$ continuity.
Integrating~\eref{eq:regularized_heaviside} gives~\cite{abbas:2010:TXF}:
\begin{equation} \label{eq:regularized_heaviside_integrated}
\psi(\phi, \varepsilon) = \frac{1}{256\varepsilon^9} \left(128\varepsilon^9 + 315\phi\varepsilon^8 - 420\phi^3\varepsilon^6 + 378\phi^5\varepsilon^4 - 180\phi^7\varepsilon^2 + 35\phi^9\right),
\end{equation}
in the region $|\phi| \le \varepsilon$. It has been suggested that a single
enrichment function may not be adequate to represent the complete range of
sharp gradients present in a problem, e.g., the shock front in a
convection-dominated problem, or the damage process zone in quasibrittle
materials. For such problems, 
multiple enrichment functions, each having a separate parameter
$\varepsilon$ should be used~\cite{patzak:2003:PZR,abbas:2010:TXF}. For the
sake of illustration consider five functions with the parameter
$\varepsilon / h = \{2.5, 0.85, 0.265, 0.085, 0.0225\}$ ($h$ is the
element size), and integrate them on the unit square ($h = 1$). This
set was produced by Abbas et al.~\cite{abbas:2010:TXF} to minimize the
pointwise error in modeling shocks. Other types of sharp-gradient enrichment
functions, such as $\tanh(q\phi)$ can also be used, where $q$ controls the
severity of the gradient. The interface is a straight line shown
in~\fref{fig:highgradient_ex1}a, and the enrichment function for
$\varepsilon = 0.085$ is depicted in~\fref{fig:highgradient_ex1}f. Four
different integration strategies are examined: (1) tensor-product over the
square, ignoring the interface; (2) and (3) integration over subtriangles
using a coarse and a fine partitioning; and (4) adaptive quadrature. In all
cases, the accuracy is improved by increasing the number of integration points
(\fref{fig:highgradient_ex1}g). The adaptive integration outperforms all other
integration strategies. The experiment is repeated for a kinked interface and
the results are shown in~\fref{fig:highgradient_ex2}. Due to the shape of the
interface, the integrand in this case is relatively more complicated, and the
advantage of using the adaptive quadrature is emphasized.
For an arbitrary curved discontinuity, different methods have been proposed,
for example, integration by partitioning the element into triangles with
curved edges~\cite{legay:2005:SAW,cheng:2010:HOX}, and
octree-subdivision~\cite{dreau:2010:SXE}. In~\fref{fig:highgradient_ex3}, the
regularized Heaviside functions are integrated over a domain with a curved
discontinuity. The interface is represented in closed form using a quadratic
polynomial. The performance of the adaptive quadrature is compared with the
tensor-product quadrature in~\fref{fig:highgradient_ex3}e. Adaptive quadratures
require fewer integration points for the same accuracy.

\begin{figure}
\setlength{\unitlength}{0.0625\linewidth}
\centering
\begin{picture}(16,10.5)
  \put(.1,7.2){\includegraphics[width=3\unitlength]{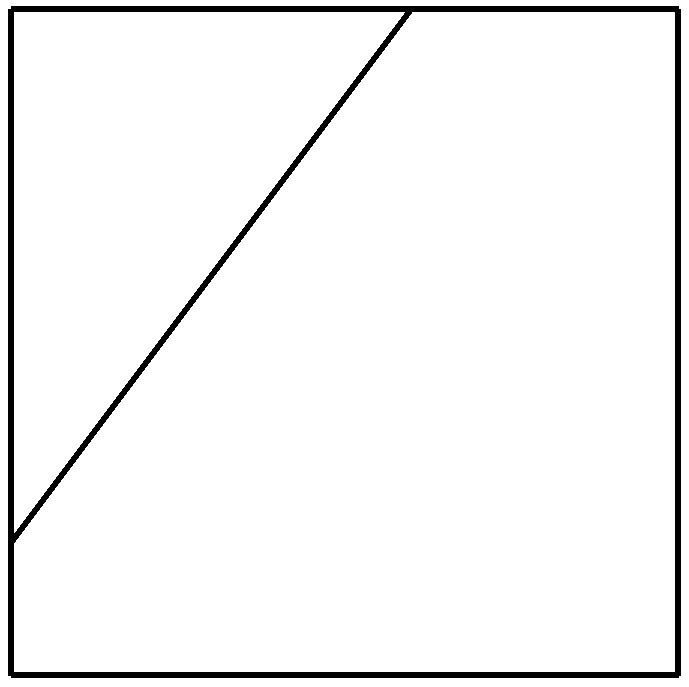}}
  \put(3.3,7.2){\includegraphics[width=3\unitlength]{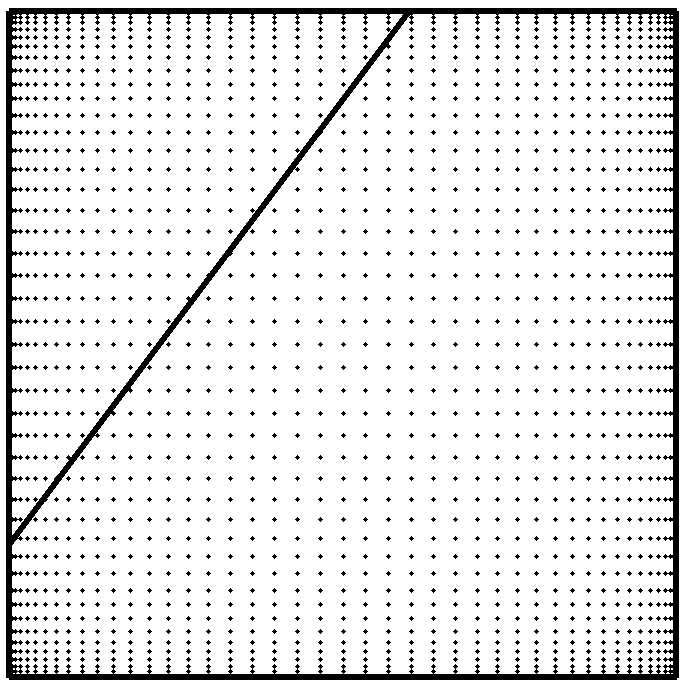}}
  \put(6.5,7.2){\includegraphics[width=3\unitlength]{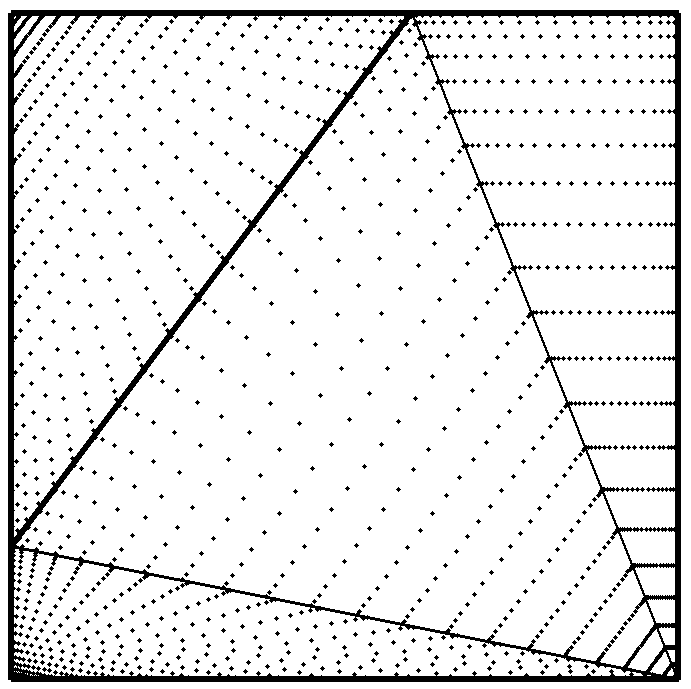}}
  \put(9.7,7.2){\includegraphics[width=3\unitlength]{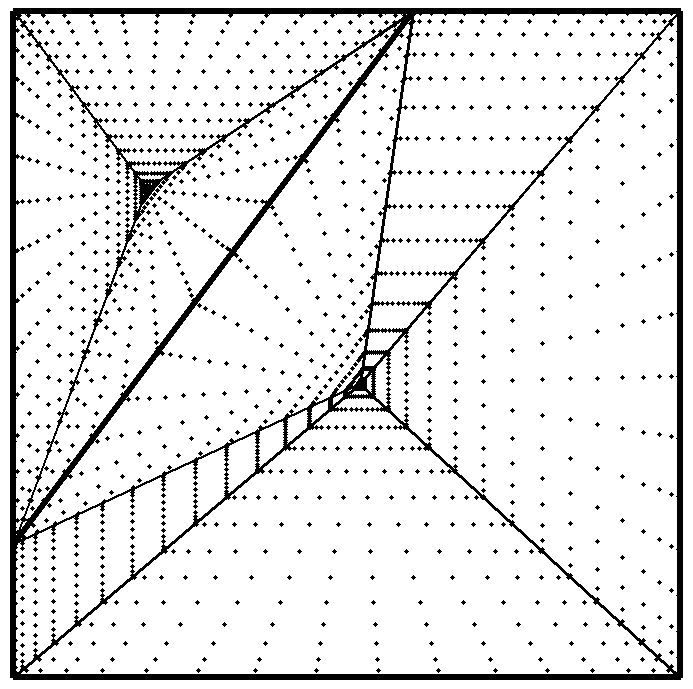}}
  \put(12.9,7.2){\includegraphics[width=3\unitlength]{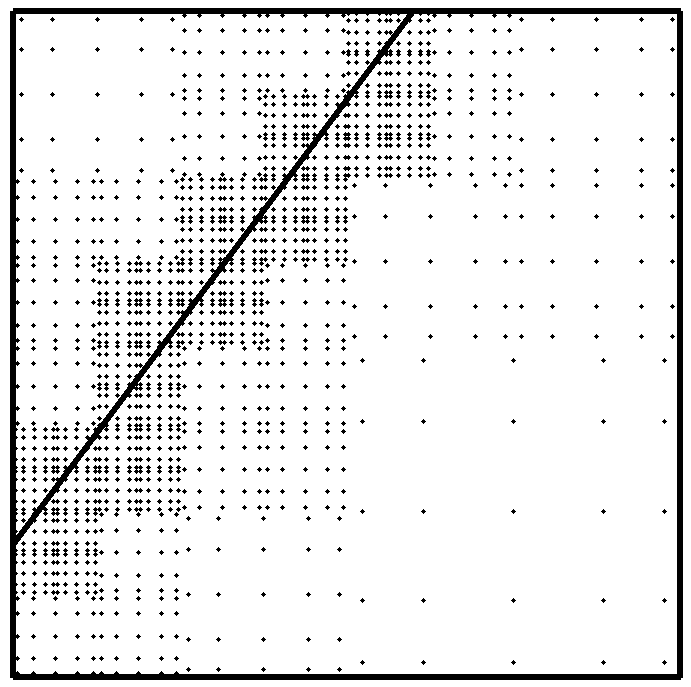}}
  \put(0,.2){\includegraphics[width=8\unitlength]{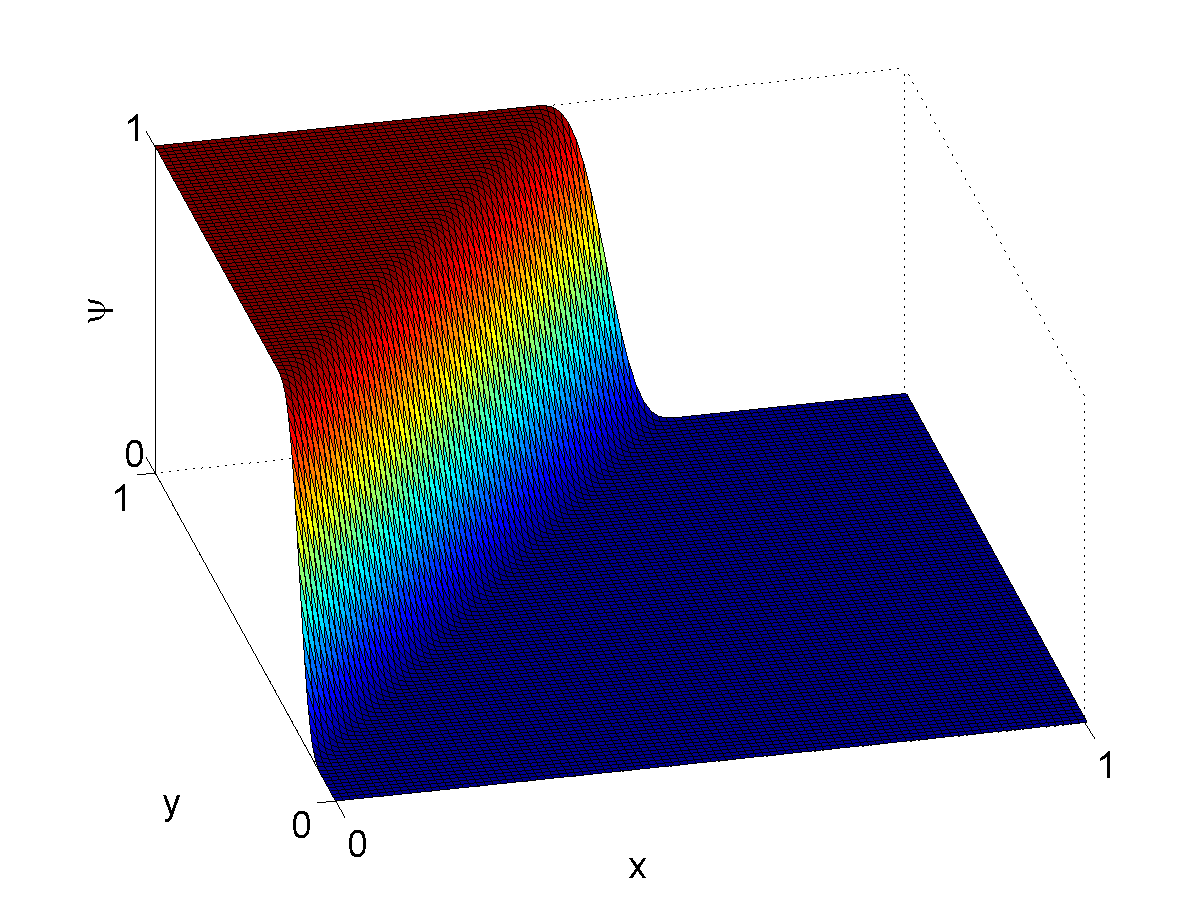}}
  \put(8,.2){\includegraphics[width=8\unitlength]{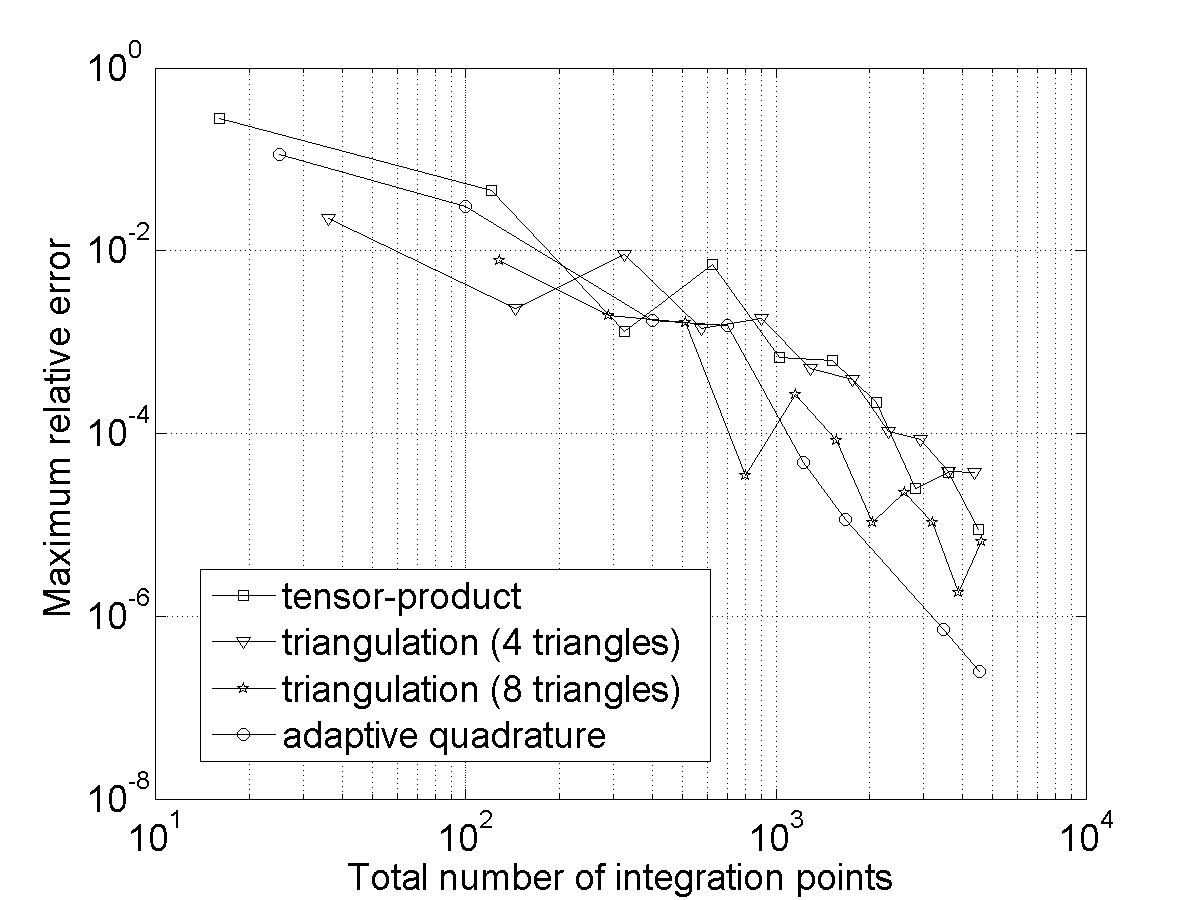}}
  \put(1.55,7){\makebox(0,0)[tc]{(a)}}
  \put(4.75,7){\makebox(0,0)[tc]{(b)}}
  \put(7.95,7){\makebox(0,0)[tc]{(c)}}
  \put(11.15,7){\makebox(0,0)[tc]{(d)}}
  \put(14.35,7){\makebox(0,0)[tc]{(e)}}
  \put(4,0){\makebox(0,0)[tc]{(f)}}
  \put(12,0){\makebox(0,0)[tc]{(g)}}
\end{picture}

\vspace*{0.1in}
\caption{Numerical integration of the regularized Heaviside function with a
straight interface.
(a) domain of integration and interface; (b) tensor-product;
(c) triangular quadratures (4 triangles);
(d) triangular quadratures (8 triangles); (e) adaptive quadrature;
(f) the regularized Heaviside function for $\varepsilon = 0.085$; and
(g) maximum relative error in the integration of the five functions versus
the number of integration points for different strategies.}\label{fig:highgradient_ex1}
\end{figure}
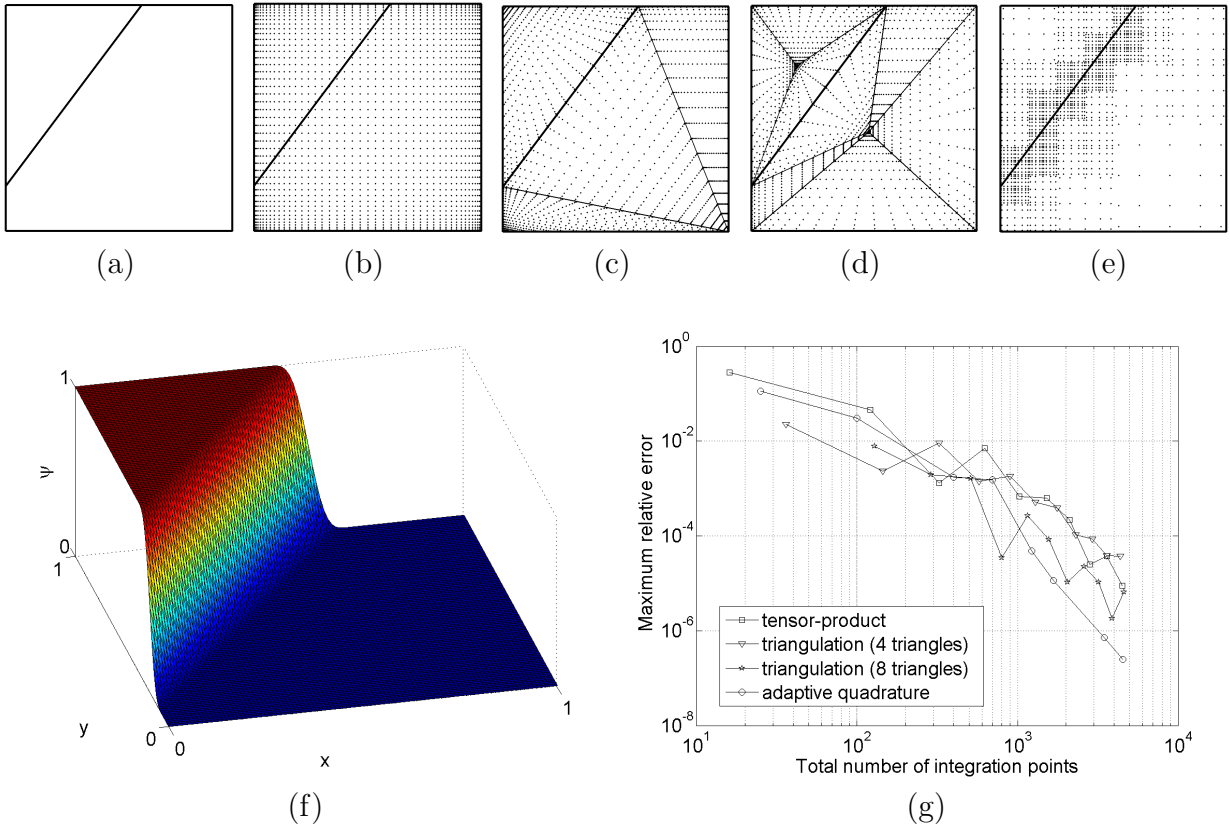
\begin{figure}
\setlength{\unitlength}{0.0625\linewidth}
\centering
\begin{picture}(16,10.5)
  \put(.1,7.2){\includegraphics[width=3\unitlength]{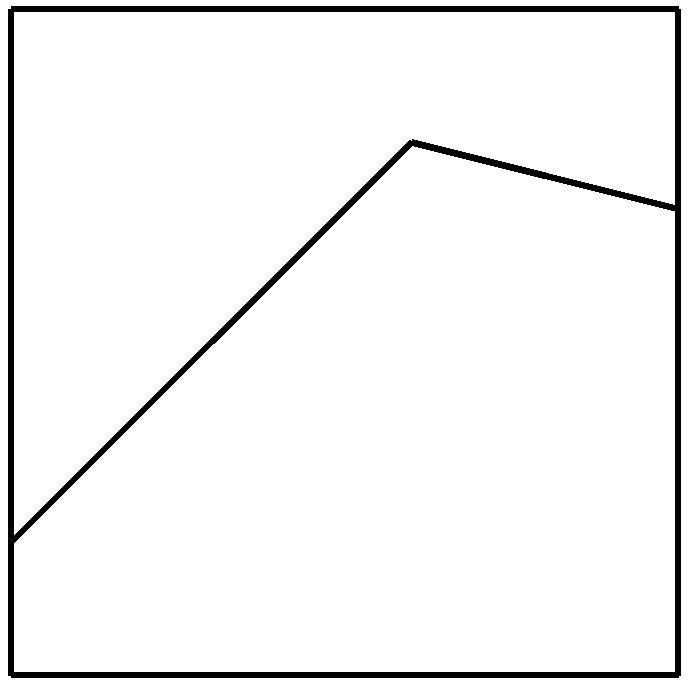}}
  \put(3.3,7.2){\includegraphics[width=3\unitlength]{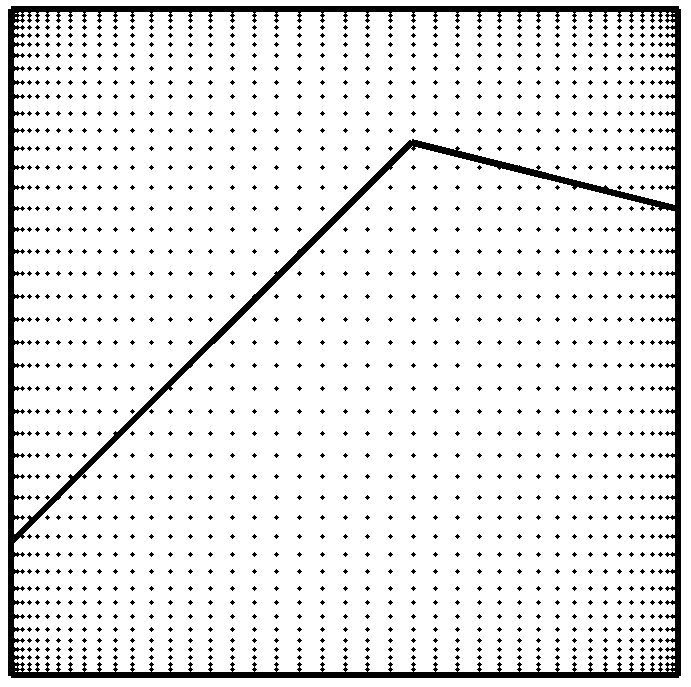}}
  \put(6.5,7.2){\includegraphics[width=3\unitlength]{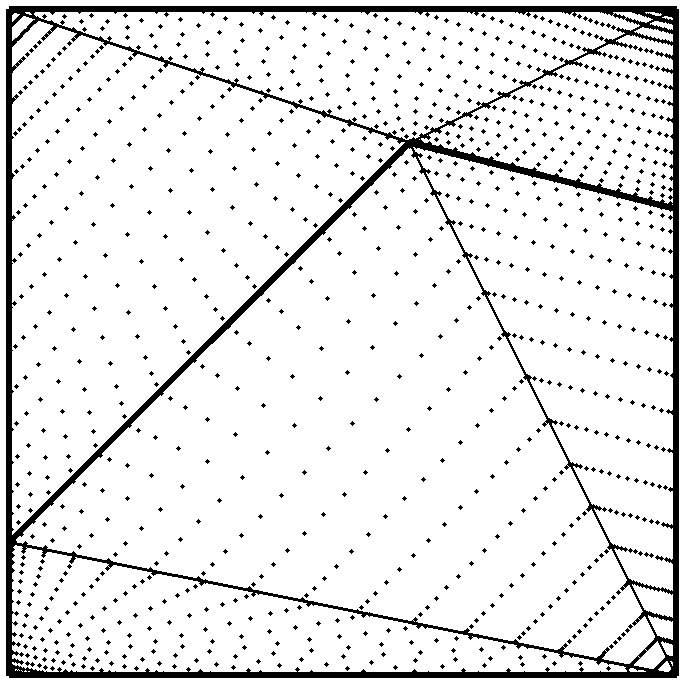}}
  \put(9.7,7.2){\includegraphics[width=3\unitlength]{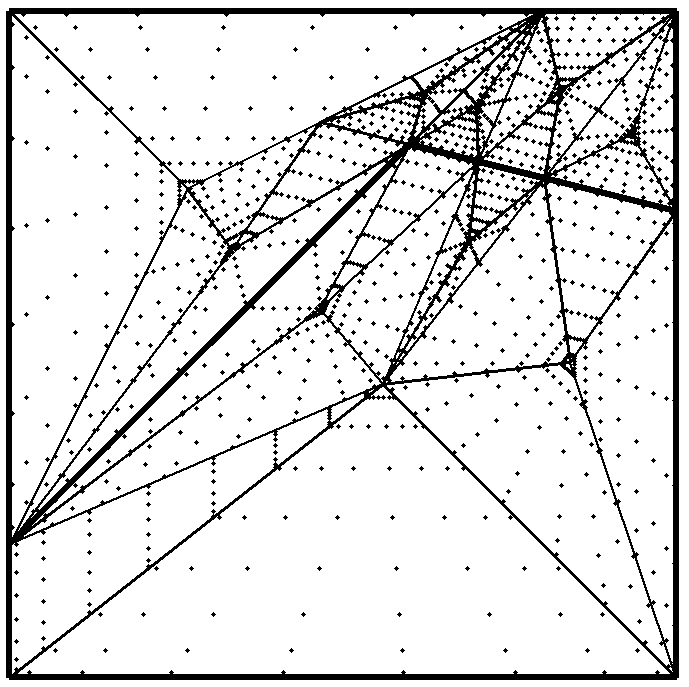}}
  \put(12.9,7.2){\includegraphics[width=3\unitlength]{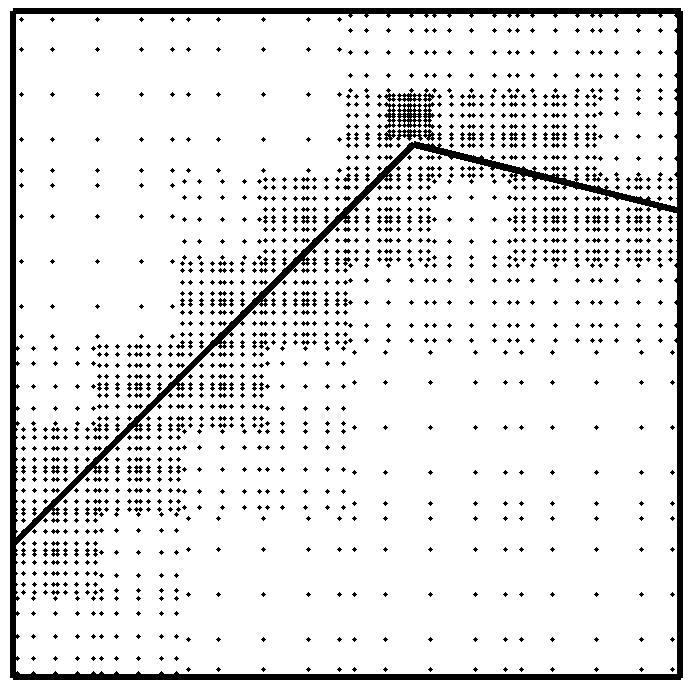}}
  \put(0,.2){\includegraphics[width=8\unitlength]{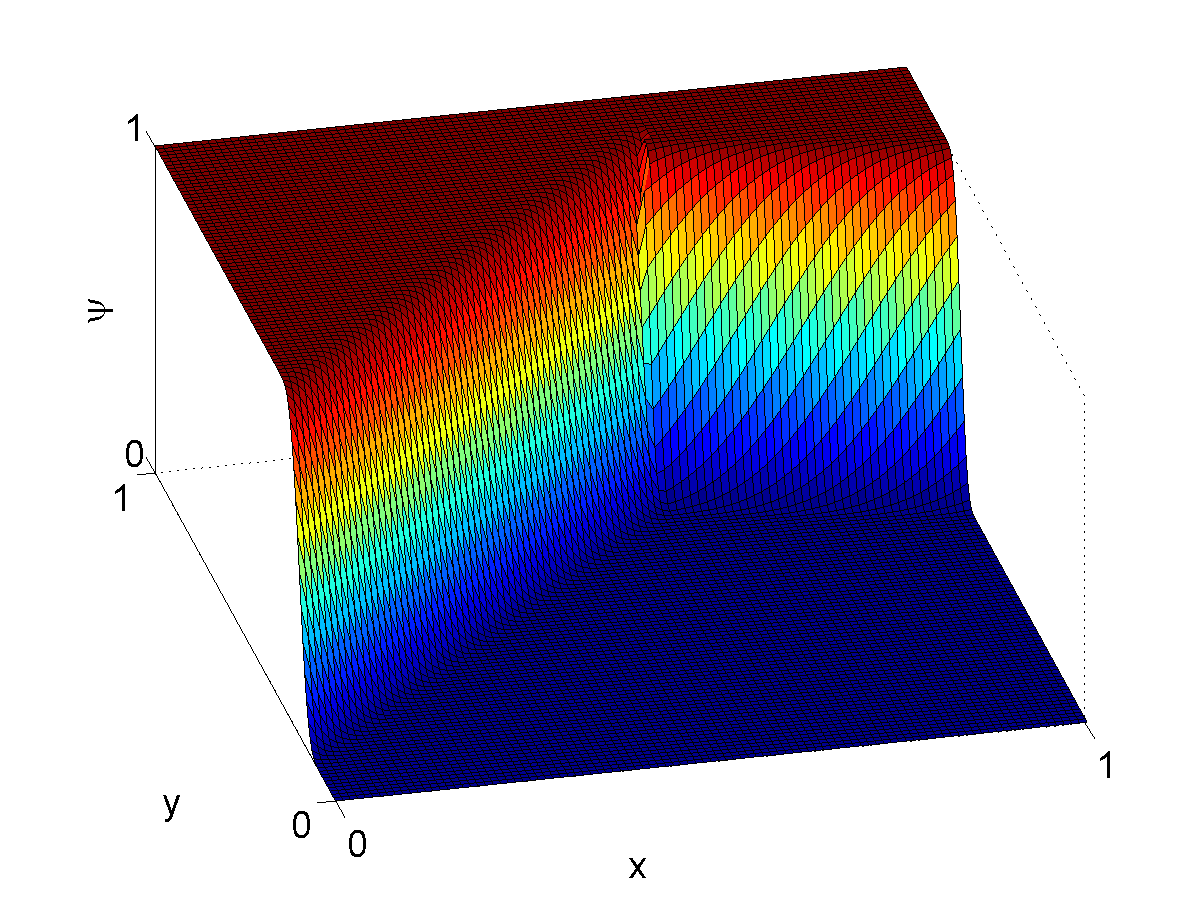}}
  \put(8,.2){\includegraphics[width=8\unitlength]{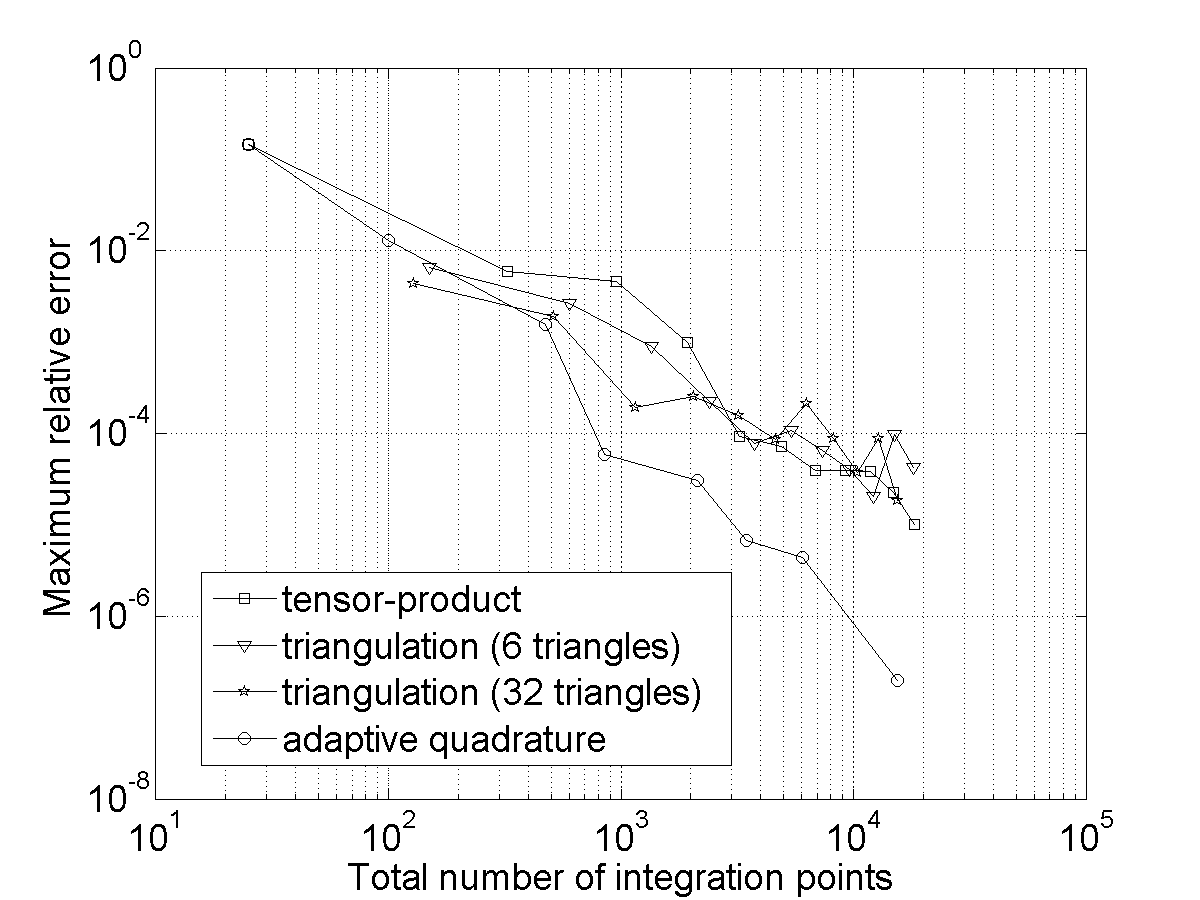}}
  \put(1.55,7){\makebox(0,0)[tc]{(a)}}
  \put(4.75,7){\makebox(0,0)[tc]{(b)}}
  \put(7.95,7){\makebox(0,0)[tc]{(c)}}
  \put(11.15,7){\makebox(0,0)[tc]{(d)}}
  \put(14.35,7){\makebox(0,0)[tc]{(e)}}
  \put(4,0){\makebox(0,0)[tc]{(f)}}
  \put(12,0){\makebox(0,0)[tc]{(g)}}
\end{picture}

\vspace*{0.1in}
\caption{Numerical integration of the regularized Heaviside function with a
kinked interface.
(a) domain of integration and interface; (b) tensor-product;
(c) triangular quadratures (6 triangles);
(d) triangular quadratures (32 triangles); (e) adaptive quadrature;
(f) the regularized Heaviside function for $\varepsilon = 0.085$; and
(g) maximum relative error in the integration of the five functions versus
the number of integration points for different strategies.}\label{fig:highgradient_ex2}
\end{figure}

\subsection{Enriched Finite Element in Quantum-Mechanical Calculations}\label{sec:poisson}
In Reference~\cite{pask:2011:LSS}, an enriched finite element method was
applied to the Coulomb problem in crystalline diamond, and tensor-product
quadratures were used for the numerical integration. The number of
integration points
was increased until convergence in the solution was achieved. A higher-order
quadrature was used over the elements containing the nuclei, since the
enriched basis functions and the source term were strongly localized about the
nuclei. A large number of integration points were required
to obtain the desired accuracy and the optimal rate of convergence. Herein,
we use the adaptive integration scheme introduced in~\sref{sec:algorithm}
to setup the system matrices, and show that the optimal rate of convergence
is realized at a relatively low computational cost.
First, we present a brief description of the problem, and then the adaptive
numerical integration algorithm is used to set up the finite element system
matrices. For details on the formulation and solution technique, see
References~\cite{sukumar:2009:CAE,pask:2011:LSS}.

\begin{figure}
\setlength{\unitlength}{0.0625\linewidth}
\centering
\begin{picture}(16,6)
  \put(0,.2){\includegraphics[width=8\unitlength]{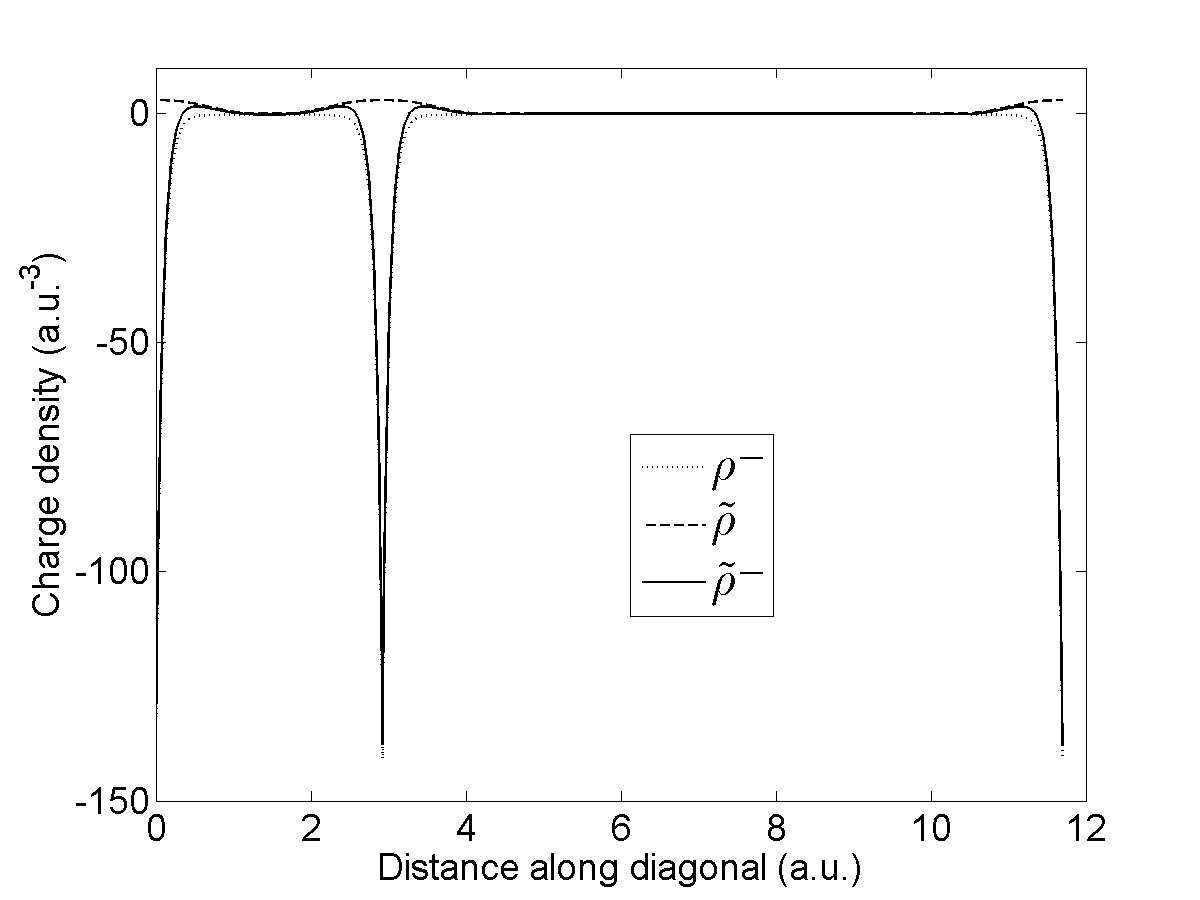}}
  \put(8,0){\includegraphics[width=8\unitlength]{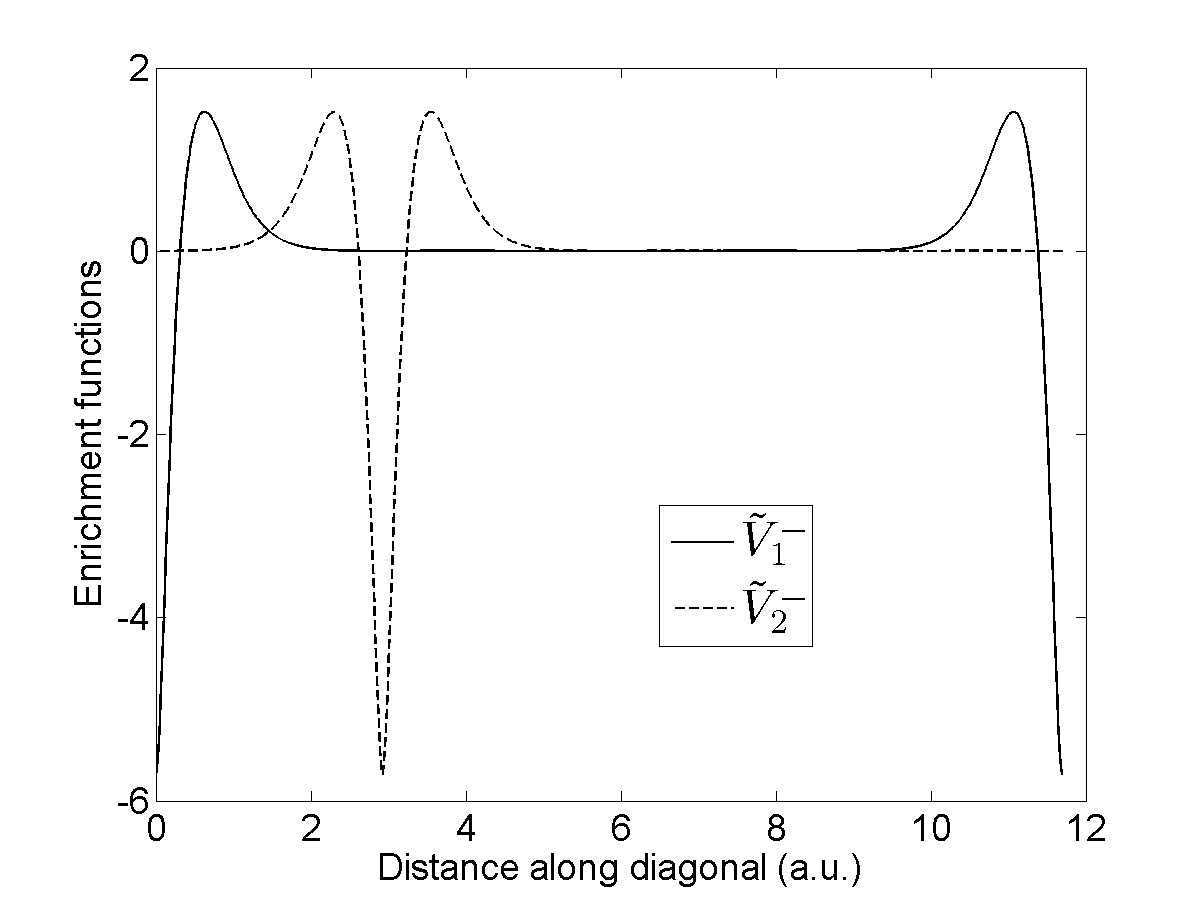}}
  \put(4,0){\makebox(0,0)[tc]{(a)}}
  \put(12,0){\makebox(0,0)[tc]{(b)}}
\end{picture}
\caption{(a) Electronic charge densities; and
(b) enrichment functions for the Coulomb problem in crystalline
diamond~\cite{pask:2011:LSS}.}\label{fig:diamond_periodic}
\end{figure}

Consider the unit cell defined by the lattice vectors
\begin{align}\label{eq:diamon_lattice}
  \vm{a}_1 &= \frac{a}{2}(0, 1, 1)  \nonumber \\
  \vm{a}_2 &= \frac{a}{2}(1, 0, 1) \\
  \vm{a}_3 &= \frac{a}{2}(1, 1, 0), \nonumber
\end{align}
where $a = 6.75$ bohr, and carbon atoms are at
$\tauv_1 = (0, 0, 0)$ and $\tauv_2 = (1/4, 1/4, 1/4)$ in lattice coordinates.
The total charge $\rho$ in the unit cell is written as
\begin{align}\label{eq:diamond_charge}
\rho(\vx) = \rhop(\vx) + \rhom(\vx)
          = \rhop(\vx) - \rhot(\vx) ~+~ \rhom(\vx) + \rhot(\vx)
          = \rhopt(\vx) + \rhomt(\vx),
\end{align}
where $\rhop(\vx) = \sum_i{\rho_i(\vx)} = \sum_i{q_i \delta(\vx-\tauv_i)}$
is the total nuclear charge density in the unit cell, $\rhom(\vx)$ is the
electronic charge density, and $\rhopt(\vx) = \rhop(\vx) - \rhot(\vx)$ and
$\rhomt(\vx) = \rhom(\vx) + \rhot(\vx)$ are the neutralized nuclear and
electronic charge densities, respectively. The neutralizing charge $\rhot(\vx)$
is introduced in~\eref{eq:diamond_charge} to circumvent the divergence of the
potential $\Vp(\vx)$ ($\Vp \sim 1/r$) at nuclear locations, so that $\Vpt$ is
extracted analytically, and $\Vmt$ is solved in real space. The
potential $\Vmt$ associated with the neutralized electronic charge density
is obtained from the solution of Poisson's equation:
\begin{equation}\label{eq:diamond_poisson}
\nabla^2 \Vmt(\vx) = - 4 \pi \rhomt(\vx),
\end{equation}
subject to periodic boundary conditions, with continuous neutralized
electronic charge density $\rhomt(\vx)$ as the source term. The EFE solution
is written as
\begin{equation}
\Vmt(\vx) = \sum_i \phi_i(\vx) a_i + \sum_\alpha{\psi_\alpha(\vx)b_\alpha}
\equiv \sum_k{\Phi_k(\vx) c_k},
\end{equation}
where $\alpha$ is summed over the atoms and
$\{\Phi_k\} = \{\phi_i\} \cup \{\psi_\alpha\}$ is the combination of the
classical and enriched basis functions that form the EFE
basis. The enrichment functions $\psi_\alpha(\vx)$ are taken as sum of the
potentials $\vImt$---isolated atomic solutions corresponding to the neutralized
electronic charge densities $\rhoImt = \rhoIm + \rhoIt$ in the vicinity of
each atom $I$. The enrichment function is written as
\begin{equation}
\psi_\alpha(\vx) = \Valmt(\vx) = \sum_\Rv \valmt(\vx-\Rv),
\end{equation}
where $\Rv$ denotes lattice translation vectors. The electronic charge
densities and the enrichment functions are plotted
in~\fref{fig:diamond_periodic}.

\begin{figure}
\setlength{\unitlength}{0.0625\linewidth}
\centering
\begin{picture}(16,7)
  \put(3,0){\includegraphics[width=10\unitlength]{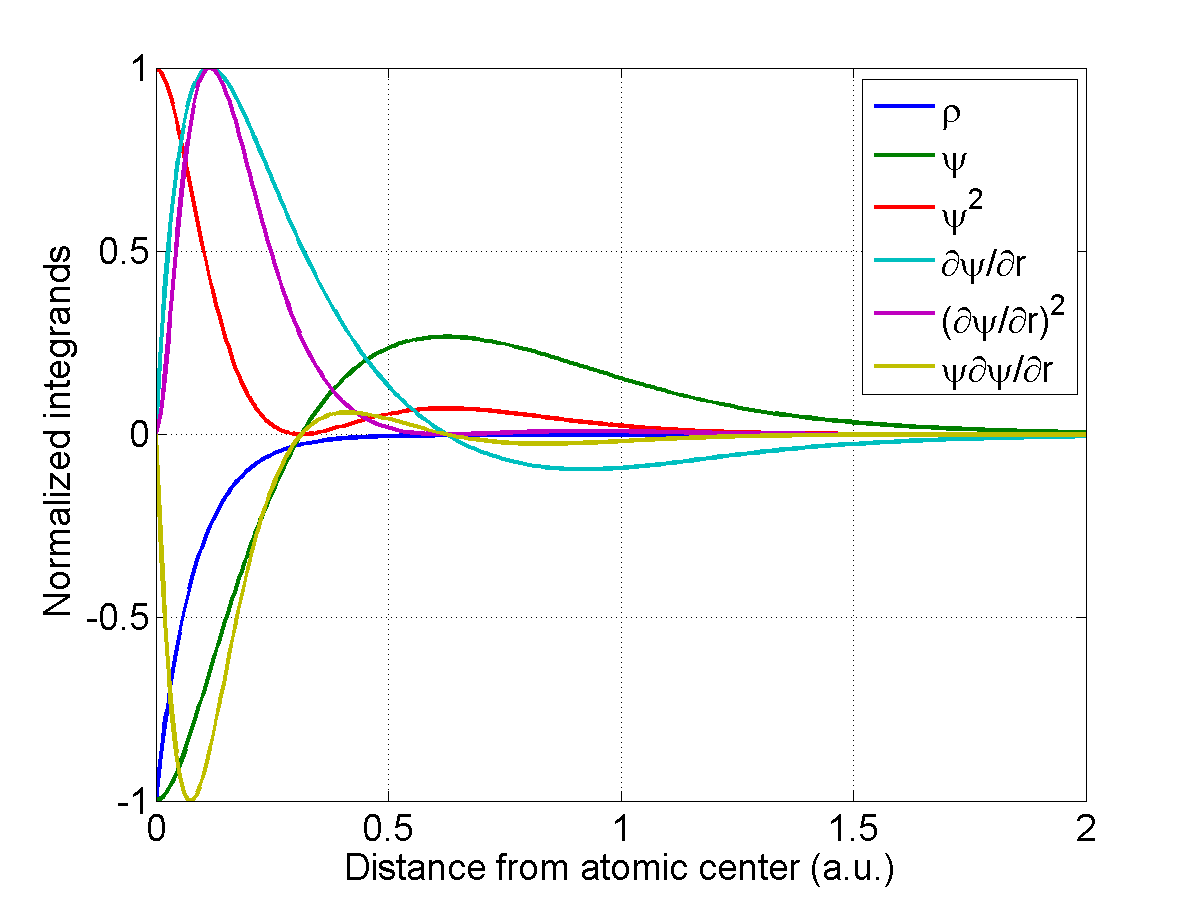}}
\end{picture}
\caption{Normalized integrands in the element stiffness matrix and force vector
of the Poisson problem.}\label{fig:diamond_all_integrands}
\end{figure}

On incorporating the trial and test functions of the
form~\eref{eq:diamond_poisson} into the weak form of the Poisson's equation,
the discrete linear system of equations emerges. It is seen that the terms
$\psi, \psi^2, \partial\psi / \partial r, (\partial\psi/\partial r)^2$, and
$\psi \partial\psi/\partial r$ appear in the element stiffness matrix, and
$\rho$ in the element force vector, where for simplicity $\psi$ and $\rho$
are used for the enrichment function, and the electronic charge density,
respectively, and $\partial\psi/\partial r$ refers to the derivative with
respect to the radial coordinate. These functions are normalized, with
an absolute maximum value of unity, and plotted
in~\fref{fig:diamond_all_integrands}.
The integrands have sharp gradients close to the atomic positions and produce
cusps at the atomic sites. While it is possible to integrate these terms
separately (resulting in multiple ad hoc quadratures), it is desirable to
have a single quadrature that is capable of efficiently evaluating the
integrals altogether.
The proposed numerical integration algorithm constructs a quadrature
rule over each finite element that satisfies a given error tolerance for all
the above integrands. Once the finite element mesh is generated, adaptive
quadratures are constructed and saved for each element, which are used
throughout the analysis and the postprocessing.

\begin{figure}
\setlength{\unitlength}{0.0625\linewidth}
\centering
\begin{picture}(16,8)
  \put(3,0){\includegraphics[width=10\unitlength]{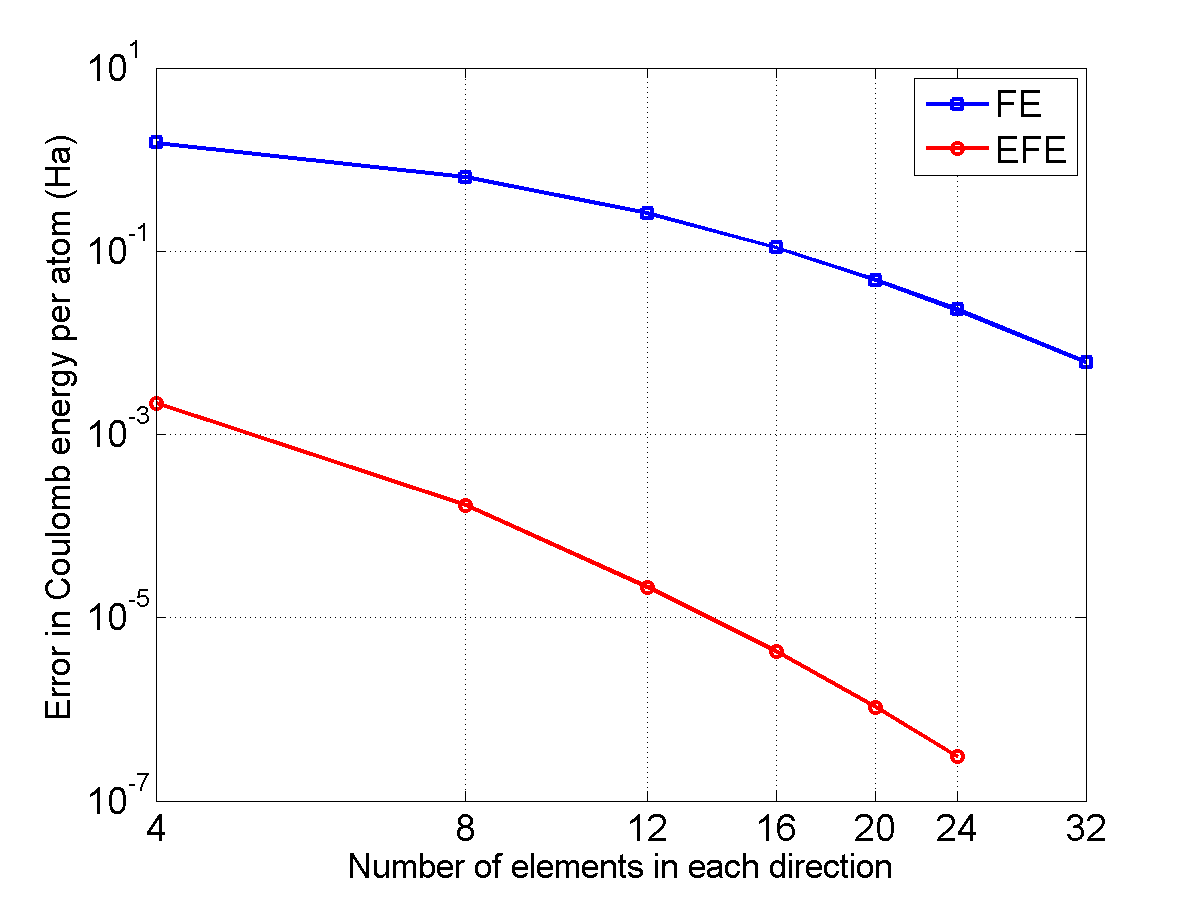}}
\end{picture}
\caption{Error in Coulomb energy per atom. $32 \times 32 \times 32$ mesh is
used as the reference solution.}\label{fig:diamond_convergence}
\end{figure}

\begin{figure}
\setlength{\unitlength}{0.0625\linewidth}
\centering
\begin{picture}(16,6)
  \put(0,.2){\includegraphics[width=8\unitlength]{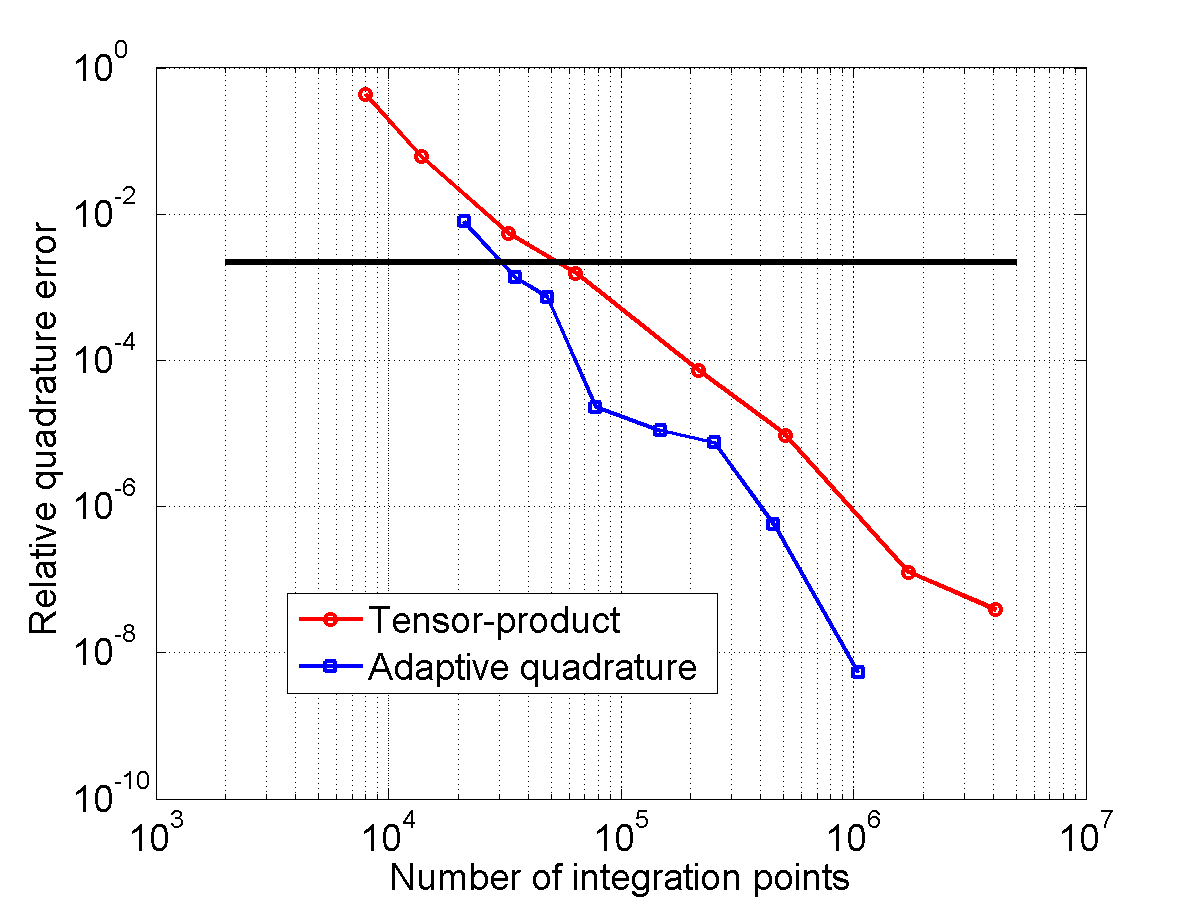}}
  \put(8,.2){\includegraphics[width=8\unitlength]{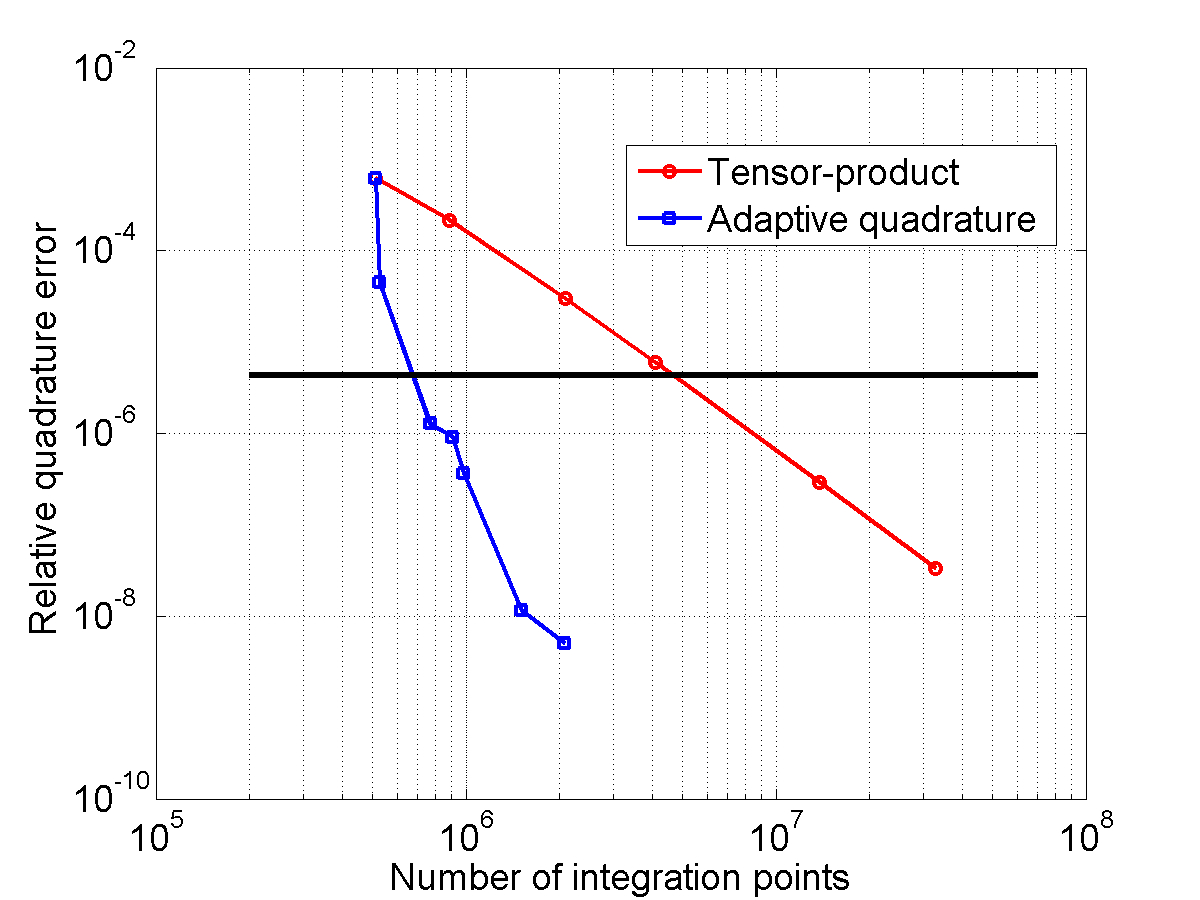}}
  \put(4,0){\makebox(0,0)[tc]{(a)}}
  \put(12,0){\makebox(0,0)[tc]{(b)}}
\end{picture}
\caption{Comparison of tensor-product and adaptive quadratures for the
crystalline diamond problem. The thick line shows the EFE solution error.
(a) $4 \times 4 \times 4$ mesh; and
(b) $16 \times 16 \times 16$ mesh.}\label{fig:diamond_comp_tensor_adap}
\end{figure}

A sequence of refined meshes are used with 4, 8, 12, 16, 20, 24, and 32 cubic
serendipity finite elements in each direction. Numerical integration is
performed with tensor-product and adaptive quadratures. In each case, the
accuracy of the quadratures is increased until convergence in the solution is
observed. The number of integration points for each mesh and integration
strategy is reported in~\tref{tab:diamond_numx}. \fref{fig:diamond_convergence}
shows the convergence curves for FE and EFE runs, with a rate of convergence
of 4.39 and 5.96, respectively, for the last three data points. The results
indicate that the integration scheme is sufficiently accurate to realize the
optimal rate of convergence.
The error tolerance (input of the quadrature construction algorithm) is the
maximum allowable error that produces a stable result (i.e., stable with
respect to further decrease in tolerance). For the pure FE runs, a $5$-point
Gauss quadrature rule is used in each direction. Adaptive integration
proves to be superior with respect to the tensor-product quadrature. The
improvement is emphasized for finer meshes where higher accuracies are
required: the integration demand of the EFE solution is only marginally
higher than the pure FE solution (of much lower accuracy) on the same mesh.
In~\fref{fig:diamond_comp_tensor_adap}, a comparison is made between the
behavior of tensor-product and adaptive quadratures for the EFE solution of the
Coulomb problem in crystalline diamond. The quadrature error in the
Coulomb energy is
plotted with respect to the number of integration points for finite element
meshes with $4$ and $16$ elements in each direction. The accuracy of the
quadrature is
increased until the desired convergence with respect to the quadrature is
attained. The adaptive quadrature requires fewer integration points to achieve
the same accuracy. Standard Gauss quadrature shows a smoother convergence
curve, which can be attributed to the uniform overall increase in the accuracy
of the quadrature throughout the domain. The efficiency of the adaptive
quadrature is more noticeable in case of the $16 \times 16 \times 16$ mesh.

Note that in all the meshes used earlier, the atoms are located at the
vertices of the finite elements---the atom inside the unit cell is at a
quarter of the diagonal from the corner, and the number of elements is a
multiple of four---and there is no atom inside any finite element. For more
general lattice systems, non-uniformly refined meshes may be required in
order to ensure that all atoms are located at the vertices of elements.
This is not desirable due to the increase in the total number of
elements and the associated computational
costs. Furthermore, with the application of the EFE, one would like to
resolve the local features by adding enrichment functions to the
approximation, and not by refining the finite element mesh. In the
following, we use
meshes with 1, 2, 3, 5, 6, and 7 elements in each direction. In all
these cases, the atom in the unit cell lies inside an element. The
number of integration points of the tensor-product and adaptive
quadratures are compared
in~\tref{tab:diamond_numx_badmesh}. The increase in the number of
integration points for the tensor-product quadrature is significant,
whereas the adaptive integration provides accurate results with a
moderate number of integration points.

\begin{table}
\caption{Comparison of the number of integration points for tensor-product
and adaptive quadratures.}\label{tab:diamond_numx}

\vspace*{0.2in}
\centering
\begin{tabular}{ll p{3.1cm} p{3.1cm} p{3.1cm}}
  \hline
  Mesh    & Error tolerance    & \multicolumn{3}{l}{Number of integration points} \\
  \cline{3-5}
          &                    & Pure FE & Tensor-product & Adaptive \\
  \hline
  4  & $2.2 \times 10^{-3}$  & 8000    & 169000      & 78000   \\
  8  & $1.7 \times 10^{-4}$  & 64000   & 624000      & 162000  \\
  12 & $2.2 \times 10^{-5}$  & 216000  & 1840000     & 349000  \\
  16 & $4.3 \times 10^{-6}$  & 512000  & 14020000    & 981000  \\
  20 & $1.1 \times 10^{-6}$  & 1000000 & 27196000    & 1549500 \\
  24 & $3.0 \times 10^{-7}$  & 1728000 & 46852000    & 3073750 \\
  32 & $4.7 \times 10^{-8}$  & 4096000 & ---         & 6112000 \\
  \hline
\end{tabular}
\end{table}
\begin{table}
\caption{Comparison of the number of integration points for tensor-product
and adaptive quadratures: one of the atoms is inside an element.}
\label{tab:diamond_numx_badmesh}

\vspace*{0.2in}
\centering
\begin{tabular}{ll p{3.1cm} p{3.1cm} p{3.1cm}}
  \hline
  Mesh    & Error tolerance    & \multicolumn{3}{l}{Number of integration points} \\
  \cline{3-5}
          &                    & Pure FE & Tensor-product & Adaptive \\
  \hline
  1  & $9.8 \times 10^{-3}$  & 125     & 1000000     & 48250   \\
  2  & $6.8 \times 10^{-3}$  & 1000    & 4096000     & 48250   \\
  3  & $3.3 \times 10^{-3}$  & 3375    & 1880000     & 55875   \\
  5  & $1.1 \times 10^{-3}$  & 15625   & 8000000     & 165250  \\
  6  & $4.9 \times 10^{-4}$  & 27000   & 27000000    & 188875  \\
  7  & $3.0 \times 10^{-4}$  & 42875   & 21952000    & 183750  \\
  \hline
\end{tabular}
\end{table}

\section{Concluding Remarks}\label{sec:conclusions}
An adaptive integration scheme was presented that can be used for
functions with sharp gradients and cusps. The algorithm uses
tensor-product quadratures with $5$ and $8$ integration points in each
direction to estimate the local error, and divides the domain uniformly,
independent of the location of the cusp or sharp gradient, until the target
tolerance is met. The error analysis in the integration of a function with
a cusp (derivative-discontinuity at a point) showed that the rate of
convergence is improved in higher dimensions.
The adaptive integration algorithm was successfully used for the integration
of a set of regularized Heaviside functions, and proved to be more efficient
than integration by partitioning as well as tensor-product quadratures.
The method was also applied to the enriched finite element solution of
the all-electron
Coulomb potential and energy of crystalline diamond (Poisson's equation).
The enrichment functions and source term were strongly localized about
the atomic positions. The adaptive integration scheme proved to be very
efficient, and recovered the optimal rate of convergence with only a
moderate increase in the number of integration points with respect to
the classical finite element method (of much lower accuracy) on the same
mesh; while reducing the integration points required in the EFE solution
by an order of magnitude or more.
The adaptive integration scheme is simple, robust, and directly
applicable to any generalized finite element method employing enrichments with
sharp local variations or cusps in $n$-dimensional parallelepiped elements.

\section*{Acknowledgements}
This work was performed in part under the auspices of the U.S. Department of
Energy by Lawrence Livermore National Laboratory under Contract 
DE-AC52-07NA27344.
S. E. Mousavi and N. Sukumar are grateful for the research support of the
National Science Foundation through contract grant DMS-0811025 to the
University of California at Davis; additional financial support from
the UC Lab Fees Research Program is also acknowledged. 

\bibliographystyle{unsrt}
\bibliography{Adap_IJNME}

\appendix
\setcounter{section}{0}
\section{\texttt{MATLAB} Code for the Adaptive Integration Scheme}\label{sec:matlab_code}
The following routine is the implementation of our adaptive integration scheme
and produces a quadrature over an $n$-dimensional hyperparallelepiped
for a given set of functions and prescribed error tolerance. First, a
description of the input and output parameters is given.
\begin{itemize}
  \item \verb|fn|: $1 \times numf$, set of integrands, defined as an array of
        structures with the member \verb|h|, which is a function handle. For
        example, in case of two integrands, we have:
        \verb|fn(1).h = @integrand1;| and \verb|fn(2).h = @integrand2;|, where
        \verb|integrand1.m| and \verb|integrand2.m| are \texttt{MATLAB}
        functions defined as $\R^n \rightarrow \R$.
  \item \verb|d|: $(n+1) \times n$, domain of integration, a
        hyperparallelepiped in $\R^n$, defined using one point as its base
        and only $n$ subsequent vertices of \verb|d|. For example, while
        \verb|d(2, :)| is a vertex of \verb|d|, the vector \verb|d(2, :)-d(1, :)|
        is a lattice vector of \verb|d|.
  \item \verb|nsp|: $1 \times 2$, number of integration points in each
        direction to evaluate the integral over the partitions. \verb|nsp(1)|
        is used to evaluate the local integrals, and \verb|nsp(2)| is used to
        evaluate the local integration error.
  \item \verb|tol|: $1 \times 1$, absolute value of the integration error.
  \item \verb|X|: $n \times numx$, quadrature points, each point is a column vector.
  \item \verb|W|: $numx \times 1$, quadrature weights.
\end{itemize}

For example, the functions used in~\sref{sec:adap_example} are defined in
the following \verb|m|-files:
{\linespread{1}
\begin{verbatim}
function val = integrand1(X)
R = sqrt(sum(X.^2, 1));
val = 10*exp(-100*R.^2);
%------------------------------------------
function val = integrand2(X)
R = sqrt(sum((X - repmat([.81, .62, .73]', 1, size(X, 2))).^2, 1));
val = 100*exp(-200*R.^2);
\end{verbatim}
}

The domain of integration (unit cube) is defined using its base and three of
its vertices:
{\linespread{1}
\begin{verbatim}
d = [0, 0, 0; ...
     1, 0, 0; ...
     0, 1, 0; ...
     0, 0, 1];
\end{verbatim}
}
With the above setting, the adaptive quadrature function can be called
(see~\fref{fig:adaptive_example}): \\
\verb|[X, W] = ndimensional_adaptive_integration(fn, d, [5, 8], 1e-6);| \\
The \texttt{MATLAB} code for the quadrature construction follows.

\linespread{1}
\begin{verbatim}
function [X, W] = ndimensional_adaptive_integration(fn, d, nsp, tol)
% External Dependencies (m-files)
% gauss_points(nsp) : 1D Gauss points
% gauss_weights(nsp): 1D Gauss weights
dim = size(d, 2); % dimension
[XYZg1, Wg1] = Get_initial_quad(nsp(1), dim);
[XYZg2, Wg2] = Get_initial_quad(nsp(2), dim);
[X, W] = Adaptive_integration(fn, d, dim, ...
             nsp, XYZg1, Wg1, XYZg2, Wg2, tol);
%------------------------------------------------------------
function [XYZ, W] = Adaptive_integration(fn, d, dim, ...
                       nsp, XYZg1, Wg1, XYZg2, Wg2, tol)
numfun = length(fn);
nspd = nsp.^dim;
PARTITION = zeros(numfun, 1);
for ind = 1:numfun
    [integ1, XYZt1, Wt1] = Integrate_over_one_cell(fn(ind), ...
                                  d, dim, nspd(1), XYZg1, Wg1);
    integ2 = Integrate_over_one_cell(fn(ind), d, dim, ...
                                    nspd(2), XYZg2, Wg2);
    err = abs(integ2 - integ1);
    if (err >= tol), PARTITION(ind) = 1; end
end
if any(PARTITION)
    two_power_dim = 2^dim;
    XYZtemp = cell(1, two_power_dim); Wtemp = XYZtemp;
    vectors = (d(2:end, :) - repmat(d(1, :), dim, 1)) / 2;
    I = zeros(dim, 1);
    for ind = 1:two_power_dim
        base = d(1, :);
        for j = 1:dim
            base = base + (I(j))*vectors(j, :);
        end
        division = [base; repmat(base, dim, 1) + vectors];
        [XYZtemp{ind}, Wtemp{ind}] = Adaptive_integration(...
                  fn(PARTITION == 1), division, dim, nsp, ...
                                XYZg1, Wg1, XYZg2, Wg2, tol);
        I = incrementI(I, dim);
    end
    XYZ = []; W = [];
    for ind = 1:two_power_dim
        XYZ = [XYZ, XYZtemp{ind}];
        W = [W; Wtemp{ind}];
    end
else
    XYZ   = XYZt1;
    W     = Wt1;
end
%------------------------------------------------------------
function [integ, XYZt, Wt] = Integrate_over_one_cell(fn, d, ...
                                           dim, nspd, XYZg, Wg)
A     = (d(2:end, :) - repmat(d(1, :), dim, 1))';
shift = d(1, :)';
XYZt  = A*XYZg + repmat(shift, 1, nspd);
Wt    = Wg*det(A);
integ = dot(Wt, fn.h(XYZt));
%------------------------------------------------------------
function [XYZg, Wg] = Get_initial_quad(nsp, dim)
% Get a nsp^dim point quad over the unit 'dim'-dimensional hypercube
xg = gauss_points(nsp)'; xg = (1 + xg)/2;
wg = gauss_weights(nsp)'; wg = wg/2;
XYZg     = xg;
Wg       = wg;
for i = 2:dim
    XYZg = repmat(XYZg, nsp, 1);
    XYZg(:, i) = reshape(repmat(xg, 1, nsp^(i-1))', 1, nsp^i);
    Wg = repmat(Wg, nsp, 1);
    Wg(:, i) = reshape(repmat(wg, 1, nsp^(i-1))', 1, nsp^i);
end
XYZg = XYZg';
Wg = prod(Wg, 2);
%------------------------------------------------------------
function I = incrementI(I, dim)
% Keep track of the vectors that construct the current cell
I(1) = I(1) + 1;
for i = 1:(dim-1)
    if (I(i) == 1), break; end
    I(i) = 0;
    I(i+1) = I(i+1) + 1;
end
\end{verbatim}

\end{document}